\pdfoutput=1
\documentclass[12pt, preprint, 3p]{elsarticle}

\usepackage{amssymb}
\usepackage{graphicx}
\usepackage{latexsym}
\usepackage{amsmath}
\usepackage{mathtools}
\usepackage{tabularx}
\usepackage{booktabs}
\usepackage{siunitx}
\usepackage[autostyle]{csquotes}
\usepackage[linesnumbered,ruled]{algorithm2e}
\newcommand*\diff{\mathop{}\!\mathrm{d}}
\usepackage{url}
\usepackage{xcolor}
\usepackage{subcaption}
\usepackage{placeins}
\usepackage{hyperref}
\usepackage{tikz, pgfplots}
\usetikzlibrary{arrows}

\journal{arXiv}

\begin{document}

\begin{frontmatter}


\title{A scale-coupled numerical method for transient close-contact melting}

\author[1,2]{Leonardo Boledi}

\author[3]{Fabian Key}
\author[1]{Benjamin Terschanski}
\author[3,4]{Stefanie Elgeti}
\author[1,2]{Julia Kowalski}

\address[1]{Chair of Methods for Model-based Development in Computational Engineering (MBD), \\ RWTH Aachen University, 52056 Aachen, Germany}
\address[2]{Aachen Institute for Advanced Study in Computational Engineering Science (AICES), \\ RWTH Aachen University, 52056 Aachen, Germany}
\address[3]{Institute of Lightweight Design and Structural Biomechanics (ILSB), \\ TU Wien, 1060 Vienna, Austria}
\address[4]{Chair for Computational Analysis of Technical Systems (CATS), \\ RWTH Aachen University, 52056 Aachen, Germany}

\begin{abstract}
We introduce a numerical workflow to model and simulate transient close-contact melting processes based on the space-time finite element method. That is, we aim at computing the velocity at which a forced heat source melts through a phase-change material. Existing approaches found in the literature consider a thermo-mechanical equilibrium in the contact melt film, which results in a constant melting velocity of the heat source. This classical approach, however, cannot account for transient effects in which the melting velocity adjusts itself to equilibrium conditions. With our contribution, we derive a model for the transient melting process of a planar heat source. We iteratively cycle between solving for the heat equation in the solid material and updating the melting velocity. The latter is computed based on the heat flux in the vicinity of the heat source. The motion of the heated body is simulated via the moving mesh strategy referred to as the Virtual Region Shear-Slip Mesh Update Method, which avoids remeshing and is particularly efficient in representing unidirectional movement. We show numerical examples to validate our methodology and present two application scenarios, a 2D planar thermal melting probe and a 2D hot-wire cutting machine.
\end{abstract}


\end{frontmatter}

\section{Introduction}
\label{sec:introduction}

The study of contact mechanics is an important topic of research for engineering applications. Many examples can be found in classical mechanical engineering, including ball bearings, gear drives and friction clutches. In recent decades, contact mechanics has been extended to new areas of development, such as fretting wear of turbine blades and material testing \cite{Popov2019}. Within the scope of this work we concentrate on close-contact melting (CCM).

CCM refers to the contact between a heated body, i.e., a heat source, and a phase-change material initially in solid state, see Fig.\ \ref{fig:ccmSetup}. Heat transfer across the contact area results in melting of the solid material and a narrow liquid-state film that develops between the heat source and the solid. The name close-contact melting underlines the physical separation between the energy source and the solid material via a micro-scale gap filled with liquid \cite{Schuller2016,Schuller2017}. Thus, heat has to be transported across the thin melt film before inducing any phase change at the liquid-solid interface. Practical applications of CCM include nuclear technology \cite{Emerman1983}, hot-wire cutting \cite{Mayer2008} and thermal heat storage systems \cite{Lacroix2001,Rozenfeld2015}. CCM theory also determines the trajectory of melting probe-type cryobots for ice exploration \cite{Zimmerman2001,Kowalski2016, Talalay2020}, which constitute one of the application scenarios considered in this study. Such thermal robots have long been employed for terrestrial research \cite{Kasser1960}, but in recent years they have been further developed in view of planetary missions \cite{Ulamec2007, Dachwald2014, Dachwald2020}. Modeling the CCM process is a key to understanding the macro-scale behaviour of the probes and to its extrapolation to extreme conditions that cannot easily be tested for in an experimental setup \cite{Schuller2017, Schuller2019}.

Available literature on CCM modeling is based on the balance laws of mass, momentum and energy in the melt film. Conservation of mass and momentum is simplified with a lubrication assumption, which results in a Reynolds equation, while the energy balance is idealized as a 2D stationary convection-diffusion equation for the temperature. The system is then closed by considering a force balance between the weight of the heat source and the film pressure, as well as a heat-flux jump condition at the solid-liquid interface that is referred to as the Stefan condition. 
For specific simplified geometries, such as spheres or ellipses, and a temperature-controlled heat source, analytical solutions for the melting velocity as a function of the exerted force can be derived \cite{Moallemi1985,Chen1994}. These results have then been extended in \cite{Schuller2016} to CCM of a temperature-controlled heat source along a curvilinear trajectory induced by a non-uniform surface temperature distribution, and to CCM of directly-powered heat sources without a dedicated temperature-control \cite{Zhao2008,Schuller2017}. The aforementioned works are, however, all based on the assumption of quasi-steady CCM, in which the melting velocity results in either a constant translational or a constant curvilinear motion. Namely, by considering an equilibrium CCM process, the energy supplied by the heat source is in balance with the sum of the latent heat required for the phase-change, the energy necessary to increase the solid material's temperature to the melting point, and the convective losses in the melt film. It is therefore not necessary to solve for the temperature profile numerically, but it can rather be subsumed as a correction factor in the so-called \emph{reduced} latent heat of melting, see for instance \cite{Schuller2016, Schuller2017}. In our contribution we want to investigate the transient CCM effect. That is, we relax the equilibrium assumption and derive a transient model, which depends on the evolving heat flux in the vicinity of the heat source. This results in a more complex scale-coupled problem, where the processes in the micro-scale melt film that determine the CCM velocity are coupled to the macro-scale temperature evolution in the surrounding phase-change material. Thus, we need to extend the existing CCM theory with a suitable numerical method to compute the transient temperature field in the solid phase-change material and to retrieve the corresponding temperature gradient, i.e., the evolving heat flux at the phase-change interface.

Finite element methods are one of the most widespread tools to numerically solve partial differential equations, e.g., the heat equation, as they can easily handle complex geometries. A particular instance of finite elements is the space-time approach, where the weak formulation is considered over a domain both in space and time \cite{Huerta2005}. This framework is especially suitable for moving boundary problems as in the CCM context, since the motion of the spatial domain over time is naturally involved \cite{Key2018}. Thus, we solve for the evolving temperature in the solid material with a time-discontinuous prismatic space-time method based on \cite{Behr1992,Behr1992.2}. Once the temperature distribution is known, the computation of the heat flux at the tip of the heat source comes into play. Most finite element schemes estimate boundary fluxes directly from the derivative of the finite element solution to the problem, which often leads to inaccuracies \cite{Akira1986}. For this reason, we employ a consistent boundary flux method, which has shown higher accuracy with respect to classical approaches \cite{Zhong1993}.

The additional challenge arises from the relative motion of the heat source with respect to the ambient phase-change material. As we are modeling a melting process, we aim at showing the displacement of the heat source over time. In order to track the moving boundaries, an update of the computational mesh is required. The simplest solution is to perform a global remeshing after the domain has been modified, which in turn significantly increases the computational effort and introduces a projection error. To avoid this error and to increase efficiency, we employ the so-called Virtual Region Shear-Slip Mesh Update Method \cite{Key2018}. Thus, we split the mesh into a static and a moving portion, where the former is related to the solid material and the latter is related to the heat source. The mesh connectivity is updated only in a small portion of the domain, which makes the representation of translational motion extremely efficient.

The article is structured as follows: In Section \ref{sec:physicalModeling}, we introduce the physical models and our workflow. The employed numerical methods, i.e., the space-time finite element formulation of the heat equation and the moving mesh method, are described in Section \ref{sec:numericalMethods}. In Section \ref{sec:numericalResults}, we present two verification cases based on the heat equation to validate our computational approach. In Section \ref{sec:2dProbe}, we show an application case of a 2D planar melting probe. In particular, the computation of the equilibrium and transient melting velocities for realistic conditions of the ambient ice is demonstrated and the temperature distribution around the probe is analyzed over time. A different application case for a hot-wire cutting process is considered in Section \ref{sec:hotWireCutting}. The results and conclusions are then summarized in Section \ref{sec:conclusions}. 

\section{Physical modeling and workflow}
\label{sec:physicalModeling}

\begin{figure}
\centering
	\includegraphics[trim=6cm 4cm 6cm 3cm,clip,width=0.7\textwidth]{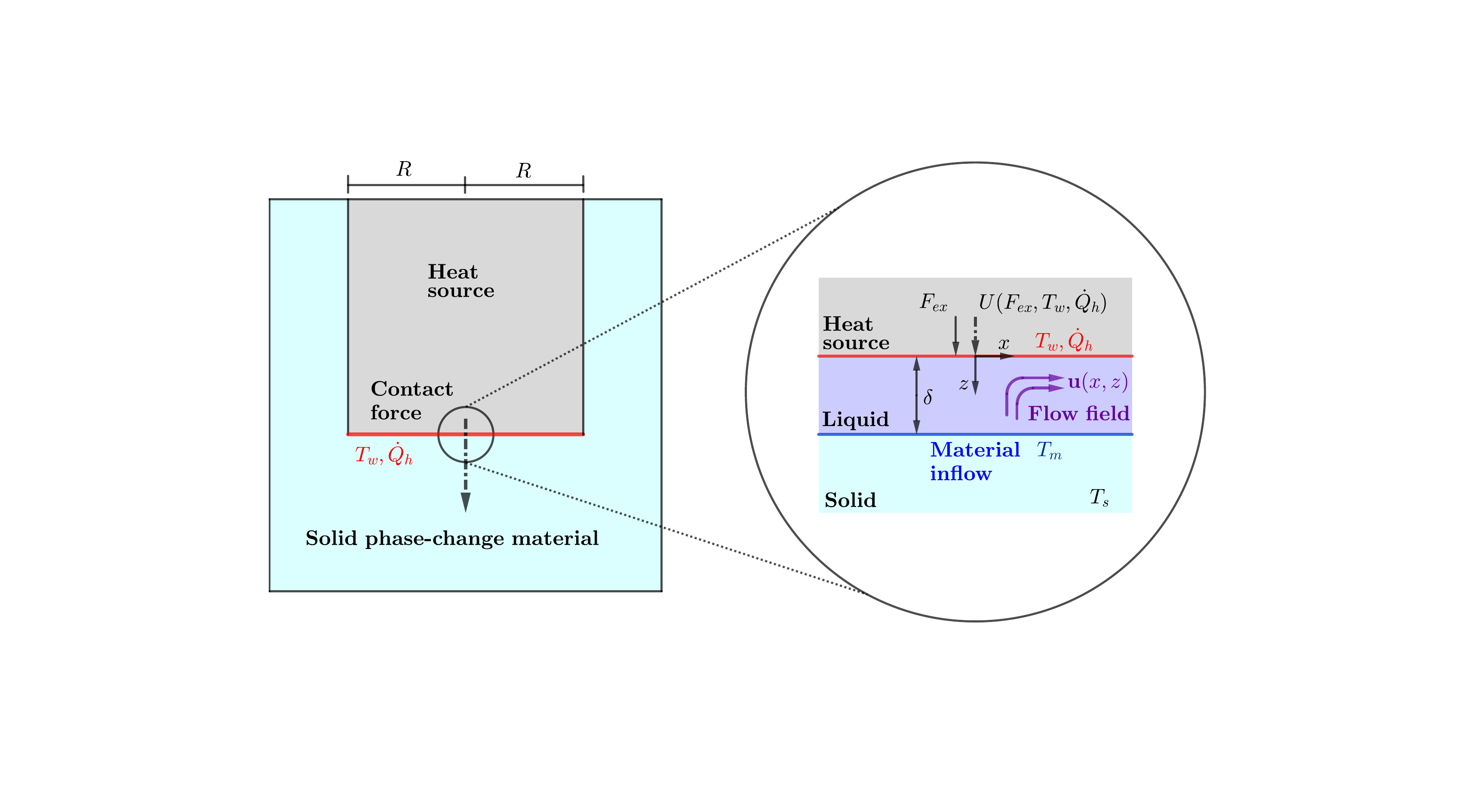}
	\caption{Sketch of the CCM process inspired by \cite{Schuller2016}. On the left, an overview of the heat source that melts through the phase-change material at velocity $U$. On the right, a close-up of the liquid melt film, note that its thickness $\delta$ is not drawn to scale. As the solid melts, liquid material is released at the phase-change interface (blue). The heat source moves downwards in response to the exerted force $F_{ex}$ and the energy released at the tip, either through temperature controlling its surface $T_w$ or by providing a heat-flow rate $\dot{Q}_h$.}
\label{fig:ccmSetup}
\end{figure}

We now introduce the mathematical models and the scale-coupled workflow for the numerical simulation of the problem. The physical situation of an idealized CCM process is sketched in Fig.\ \ref{fig:ccmSetup}. The planar heat source is separated from the solid phase-change material via the liquid melt film, which is sustained by a balance between exerted weight and liquid film pressure. Note that the melt-film thickness $\delta$ is not depicted to scale. The heat source releases energy into the liquid film, which leads to melting of the solid and subsequently leads to a downwards motion of the heat source indicated by $U$. The heat source can either be temperature-controlled, such that the heat source's surface temperature $T_w$ is given, or its heat-flow rate $\dot{Q}_h$ is prescribed. 

\begin{figure}
\centering
	\includegraphics[trim=0 0.0cm 0 0.0cm,clip,width=0.9\textwidth]{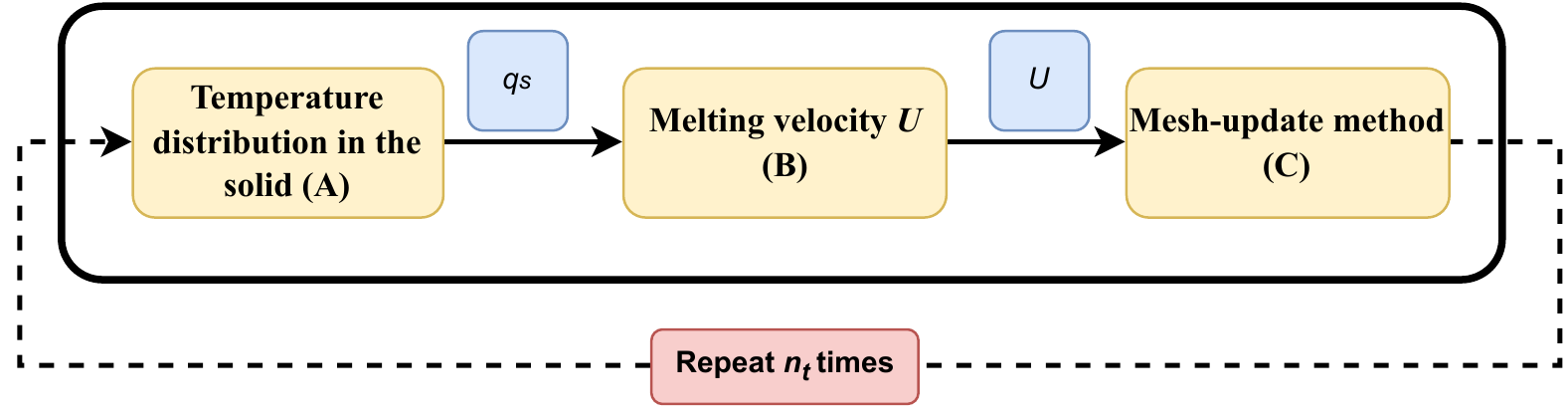}
	\caption{Overview of the scale-coupled computational approach at each time step. The macro-scale temperature distribution \textbf{(A)} around the heat source is computed via space-time finite elements. The melting velocity \textbf{(B)} is then retrieved according to the micro-scale model in Sections \ref{sec:eqVelocity} and \ref{sec:transientVelocity}, where the latter requires the numerical heat flux $q_s$ at the phase interface. Finally we apply the velocity $U$ to the mesh movement \textbf{(C)}.}
\label{fig:workflow}
\end{figure}

In order to solve for the CCM problem, we need to retrieve the evolving temperature profile around the moving heat source and the corresponding melting velocity. Fig.\ \ref{fig:workflow} gives an overview of the essential building blocks of our solution approach at each time step. We underline three main steps: \textbf{(A)} modeling the evolution of temperature in the ambient phase-change material, \textbf{(B)} the computation of the melting velocity $U$ based on an embedded analytical model accounting for processes in the micro-scale melt film and \textbf{(C)} the mesh update method to account for the downwards movement of the heat source.

Starting with block \textbf{(A)}, we consider an Eulerian frame of reference to model the spatio-temporal temperature field around the heat source. That is, we represent the solid as a static material through which the heated body moves. Let $\Omega$ be the computational domain and $(0,t_{\text{end}})$ a time interval. The evolution of the temperature field $T(\textbf{x},t)$ is then governed by the unsteady heat equation
\begin{equation}
    \frac{\partial T}{\partial t}-\frac{\kappa_s}{\rho_s \, c_{p,s}}\,\Delta T = 0 \hspace{5mm} \text{in} \hspace{2mm} \Omega\times(0,t_{\text{end}}),
    \label{eq:heat}
\end{equation}
with thermal conductivity $\kappa_s$, density $\rho_s$ and specific heat capacity $c_{p,s}$. To increase readability, we omit the dependency of each term on space and time $(\textbf{x},t)$. Initial and boundary conditions are necessary to solve the system. In our case, we will distinguish between the boundary given by the heat source surface, at which an analytical model is embedded to account for the  processes in the micro-scale melt film, see following subsection, and standard Dirichlet- or Neumann-type boundaries of the computational domain. The numerical approach to solve Eq.\ \eqref{eq:heat} is described in Section \ref{simulationTemperature}. In the following, we introduce the analytical models to consider the processes in the micro-scale melt film and estimate the melting velocity of the heat source.

\subsection{Close-contact melting velocity at equilibrium conditions}
\label{sec:eqVelocity}

The melting velocity, i.e., the velocity at which the heat source moves downwards into the phase-change material, has two roles in our model: First, it determines the motion of the heated body, hence it is a relevant quantity to update the computational domain; Second, it determines the mass intake rate into the micro-scale melt film. In order to to derive the melting velocity of the heat source, that is block \textbf{(B)} in Fig.\ \ref{fig:workflow}, we analytically model heat and flow fields in the melt film with external CCM theory, see Fig.\ \ref{fig:ccmSetup}. Due to the weight that the heated body exerts on the melt film, the latter will be very thin and present itself as a lubricated layer. CCM theory therefore combines lubrication theory with a Stefan-type phase-change model. Assuming idealized heat source geometries \cite{Chen1994}, the problem can be analytically solved for the melting velocity in response to exerted weight and temperature/heat flux forcing at the heat source.

Within the scope of this work, we consider a 2D planar heat source of either constant temperature or constant heat-flow rate at its surface, which results in a unidirectional motion of the heated body. Here, we briefly sketch the analytical model evaluated at the heat source boundary, which follows \cite{Hamrock2004}. The reference coordinates for the derivation, fixed to the heat source, are depicted in Fig.\ \ref{fig:ccmSetup}. The model is based on the assumption that the melt-film thickness $\delta$ is much smaller than the characteristic length $R$ of the heat source and that inertia terms are negligible, i.e., 
\begin{equation}
    \varepsilon=\frac{\delta}{R}\ll 1, \hspace{3mm} Re\left(\frac{\delta}{R}\right)^2\ll 1.
    \label{eq:lubTheory}
\end{equation}
The Reynolds number $Re$ can be expressed in terms of the characteristic velocity $u_0$ in the $x$-direction tangential to the surface of the heat source
\begin{equation}
    Re=\frac{\rho_l \, u_0 \, R}{\mu_l}.
    \label{eq:reynolds}
\end{equation}
Here, $\rho_l$ and $\mu_l$ denote the density and viscosity of the liquid, respectively. We now list the most important steps in deriving the melting velocity according to CCM theory. We refer to \cite{Schuller2016, Schuller2019} for additional details:
\begin{enumerate}
    \item Formulate balance laws for mass, momentum and energy in the melt film with respect to a Lagrangian reference frame that is fixed to the heat source with velocity $U$, cf.\ Fig.\ \ref{fig:ccmSetup}. Use lubrication Eqs.\ \eqref{eq:lubTheory} and \eqref{eq:reynolds} to simplify the governing equations;
    \item Assume a no-slip velocity boundary condition at the surface of the heat source. At the phase-change interface we have mass inflow due to the ongoing melting process. That is, a Dirichlet condition for the flow field is set, which is proportional to the melting velocity $U$. The melting temperature $T_m$ is then imposed the phase interface. Finally, The heat source boundary condition is chosen according to either a temperature-controlled or heat flux-controlled setting;
    \item Integrate the momentum equation across the melt film to obtain the horizontal velocity as a function of the pressure gradient and melt-film thickness $\delta$. Substitution of the boundary conditions yields an expression for the pressure gradient and the horizontal velocity as a function of $\delta$;
    \item Compute the melt-film thickness $\delta$ by assuming a parametrized temperature profile, for instance linear \cite{Schuller2016} or quadratic \cite{Schuller2019}, across the melt film and by considering local energy balance across the phase interface known as the Stefan condition \cite{Stefan1891}, that is
    \begin{equation}
    \kappa_l\,\frac{\partial T}{\partial z}\bigg|_{z=\delta}=-\rho\, U\left[h_m+c_{p,s}\,(T_m-T_s)\right]=-\rho\,U\,h_m^*.
    \label{eq:eqStefan}
\end{equation}
    In the expression $h_m$ denotes the latent heat of melting, $\kappa_l$ denotes the thermal conductivity and $c_{p,s}$ denotes the heat capacity. The two terms $T_m$, $T_s$ represent the melting temperature and the initial temperature of the solid material. The subscripts $(\cdot)_{l,s}$ denote material properties associated with the liquid phase and the solid phase, respectively. Note that the reduced latent heat of melting $h_m^*$ already accounts for the necessary energy to bring the solid material to the melting temperature, so that we do not need to numerically solve for the temperature field yet;
    \item The melting velocity $U$ can now be determined via a force balance in the melt film, hence balancing the weight exerted by the heat source with the melt-film pressure. 
\end{enumerate}
Note that we assume a constant melt-film thickness in tangential $x$-direction following \cite{Schuller2016, Schuller2019}. This assumption is appropriate for the considered planar heat source, yet it can be relaxed in principle, see for instance \cite{Schuller2017}.

Following the aforementioned steps, we derive the melting velocity $U$, which under equilibrium conditions is constant over time.
If we impose a constant temperature value $T_w$ on the heat source surface, we end up with the relation 
\begin{equation}
    U=\left(\frac{[(T_w-T_m)\,\kappa_l]^3}{8\mu_l\,(\rho_s\,h_m^*\,R)^3} \, F_{ex}\right)^\frac{1}{4}
    \label{eq:meltTemp}
\end{equation}
that we call \textit{temperature-controlled} CCM velocity \cite{Schuller2016}. Alternatively, the heat-flow rate $\dot{Q}_h$ at the heat source can be prescribed, which yields
\begin{align}
\begin{split}
    &\frac{\rho_s\,U\,h_m^*}{q_h}\left(\frac{7F(U)}{20\alpha_l}+1\right)+\frac{3F(U)}{20\alpha_l}-1=0, \\ &F(U)=\left(\frac{\rho_s}{\rho_l}\,R\,U\right)^\frac{4}{3}\left(\frac{3\pi\,\mu_l}{2F_{ex}}\right)^\frac{1}{3},
\end{split}
\label{eq:meltPower}
\end{align}
i.e., the \textit{power-controlled} CCM velocity \cite{Schuller2019}. Here $q_h$ is the bulk heat-flow rate of the heat source $\dot{Q}_h$ per unit area. The power-controlled CCM velocity is given as a non-linear relation for the velocity $U$, which we solve with a basic secant method. Note that both CCM velocity Eqs.\ \eqref{eq:meltTemp} and \eqref{eq:meltPower} depend on the force $F_{ex}$ exerted by the heat source onto the melt film. We refer to \cite{Schuller2016} for a detailed discussion of such term, but we underline that its value is known.
The takeaway here is that the melting velocity, which will later determine the heat source's solid body motion in the grid, can be computed as a function of the temperature imposed on the heat source in the temperature-controlled case, or as a function of the prescribed heat-flow rate in the power-controlled case.

\subsection{Transient close-contact melting}
\label{sec:transientVelocity}

Past studies allow to determine the equilibrium CCM velocity at which the heat source sinks into the phase-change material, yet they cannot account for the expected transient velocity ramp-up in response to a sudden or gradual change of the heat source's temperature or heating power. To include these transient effects, we relax the equilibrium assumption of Section \ref{sec:eqVelocity} and consider the full Stefan condition \cite{Stefan1891} at the interface between the liquid melt film and the solid material, that is
\begin{equation}
    \kappa_l\,\frac{\partial T}{\partial z}\bigg|_{\delta^-}-\underbrace{\kappa_s\,\frac{\partial T}{\partial z}\bigg|_{\delta^+}}_{-q_s}= -\rho_s\,h_m\,U.
    \label{eq:stefanGen}
\end{equation}
The term $q_s$ denotes the time-dependent heat flux at the solid side of the phase interface. This quantity can be extracted from the macro-scale temperature evolution in the ambient phase-change material, see Eq.\ \eqref{eq:heat}. We will address its numerical approximation in Section \ref{sec:cbfMethod}. 

Assuming that we know the heat flux at the solid side of the interface, we can modify Eqs.\ \eqref{eq:meltTemp} and \eqref{eq:meltPower} according to Eq.\ \eqref{eq:stefanGen} and obtain
\begin{equation}
    F_{ex} - \frac{8\mu_l\,U(\rho_s\,U\,h_m+q_s)^3R^3}{[(T_w-T_m)\,\kappa_l]^3} = 0
    \label{eq:meltTempTransient}
\end{equation}
for the temperature-controlled setup and 
\begin{align}
\begin{split}
    &\left(\frac{\rho_s\,h_m\,U+q_s}{q_h}\right)\cdot\left(\frac{7F(U)}{20\alpha_l}+1\right)+\frac{3F(U)}{20\alpha_l}-1 = 0, \\ &F(U)=\left(\frac{\rho_s}{\rho_l}\,R\,U\right)^\frac{4}{3}\left(\frac{3\pi\,\mu_l}{2F_{ex}}\right)^\frac{1}{3}
\end{split}
\label{eq:meltPowerTransient}
\end{align}
for the power-controlled setup. By inserting the equilibrium flux from Eq.\ \eqref{eq:eqStefan} in Eq.\ \eqref{eq:meltPowerTransient}, we immediately see that this result is consistent with Eq.\ \eqref{eq:meltPower} for the equilibrium melting velocity.

The analytical relation of the melting velocity now depends on the temperature gradient in the ambient phase-change material. It can hence be integrated into a scale-coupled Eulerian-Lagrangian approach, see again Fig.\ \ref{fig:workflow}. Namely, we first compute the transient temperature field in the ambient phase-change material in an Eulerian reference frame, then we estimate the melting velocity by extracting the heat flux at the heat source surface. This velocity is analytically modelled by considering the micro-scale processes in the melt film as sketched before. The micro-scale melt-film model follows CCM theory and is formulated in a Lagrangian reference frame. In the third step, the computational domain is updated according to the derived mesh velocity.

\section{Numerical methods}
\label{sec:numericalMethods}

After having presented the physical foundations to model the motion of the heat source through a solid phase-change material, we now discuss the computational approach that we take for that purpose. The movement of the heated body and the underlying temperature distribution in the solid are determined by a numerical approximation of the transient governing equations presented above, making use of the velocity modeling shown in the previous section. In particular, we perform an iterative procedure based on three stages:
\begin{enumerate}
    \item Compute the development of the temperature field due to the heat source;
    \item Estimate the resulting velocity of the melting process;
    \item Incorporate this movement in our simulation.
\end{enumerate}
Here, each iteration of the procedure refers to a time interval from the overall period of observation. In the following, each individual step will be described in more detail.

\subsection{Simulation of the temperature field}
\label{simulationTemperature}

In order to determine the transient temperature distribution in the solid phase according to Eq.\ \eqref{eq:heat}, we employ a space-time finite element method (ST-FEM). We select ST-FEM in this application, since both the temperature field itself and the computational domain evolve over time, the latter due to the movement of the heat source through the domain.  
\par
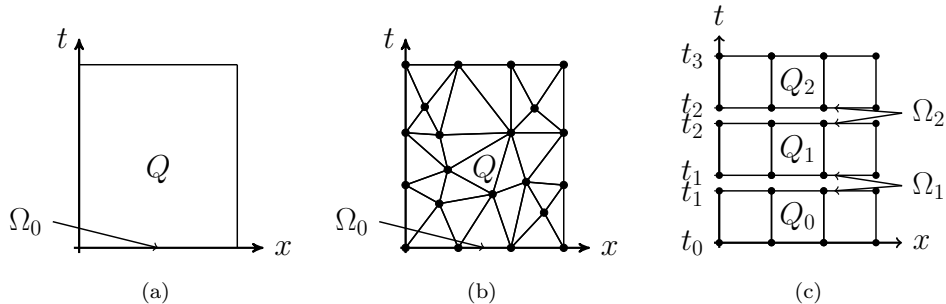
\begin{figure}
    \Large
    \centering
    \begin{subfigure}[b]{0.25\textwidth}
        \resizebox{\textwidth}{!}{\begin{tikzpicture}[
    axis/.style={very thick, ->, >=stealth'},
    important line/.style={thick},
    pile/.style={thick, ->, >=stealth', shorten <=2pt, shorten
    >=2pt},
    every node/.style={color=black}
    ]
    \draw[axis] (-0.1,0)  -- (3.5,0) node(xline)[right]{$x$};

    \draw[axis] (0,-0.1) -- (0,4) node(yline)[left]{$t$};
  	\node (Q) at (1.5,1.5) {$Q$};
    \draw[thick,<-] (1.5,0cm) -- (-0.5, 0.5) node[anchor=east] {$\Omega_0$};
    \draw[thick] (3.0, 0.0) -- (3.0, 3.5);
    \draw[thick] (0.0, 3.5) -- (3.0, 3.5);
\end{tikzpicture}}
        \subcaption{}
        \label{fig:spaceTimeDomain}
    \end{subfigure}
    \begin{subfigure}[b]{0.25\textwidth}
        \resizebox{\textwidth}{!}{\begin{tikzpicture}[
    axis/.style={very thick, ->, >=stealth'},
    important line/.style={thick},
    pile/.style={thick, ->, >=stealth', shorten <=2pt, shorten
    >=2pt},
    every node/.style={color=black}
    ]
    \draw[axis] (-0.1,0)  -- (3.5,0) node(xline)[right]{$x$};

    \draw[axis] (0,-0.1) -- (0,4) node(yline)[left]{$t$};
\begin{axis}[axis lines=none,
	xmin = 0,
	xmax = 3,
	ymin =0,
	ymax = 3.5,
	x=1cm,
	y=1cm
	]
	\addplot [patch, patch refines=0, mesh, thick, black, patch table={images/d1_triangles.dat}] table {images/d1_points.dat};
  \addplot [mark=*, draw=black, fill=black, only marks] table {images/d1_points.dat};
\end{axis}
  \node (Q) at (1.5,1.5) {$Q$};
  \draw[thick,<-] (1.5,0cm) -- (-0.5, 0.5) node[anchor=east] {$\Omega_0$};
  \draw[thick] (3.0, 0.0) -- (3.0, 3.5);
  \draw[thick] (0.0, 3.5) -- (3.0, 3.5);
\end{tikzpicture}}
        \subcaption{}
        \label{fig:spaceTimeDomainCSST}
    \end{subfigure}
    \begin{subfigure}[b]{0.25\textwidth}
        \resizebox{\textwidth}{!}{\newcommand*\cols{2}
\newcommand*\rows{3}
\begin{tikzpicture}[
    scale=1,
    axis/.style={very thick, ->},
    important line/.style={thick},
    every node/.style={color=black}
    ]
    \draw[axis] (-0.1,0)  -- (3.5,0) node(xline)[right]{$x$};
    \draw[very thick] (0,-0.1) -- (0,1);
    \draw[very thick] (0,1.3) -- (0,2.3);
    \draw[axis] (0,2.6) -- (0,4) node(yline)[above]{$t$};
    \draw (0.1,0cm) -- (-0.1, 0cm) node[anchor=east] {$t_{0}$};
    \draw (0.1,1cm) -- (-0.1, 1cm) node[anchor=east, yshift=-0.05cm] {$t_{1}$};
    \draw (0.1,1.3cm) -- (-0.1, 1.3cm) node[anchor=east,yshift=0.05cm] {$t_{1}$};
    \draw (0.1,2.3cm) -- (-0.1, 2.3cm) node[anchor=east, yshift=-0.05cm] {$t_{2}$};
    \draw (0.1,2.6cm) -- (-0.1, 2.6cm) node[anchor=east, yshift=0.05cm] {$t_{2}$};
    \draw (0.1,3.6cm) -- (-0.1, 3.6cm) node[anchor=east] {$t_{3}$};
    \foreach \row in {0, 1.3, ...,\rows} {
     		\foreach \col in {0, 1, ...,\cols} {
   				 \draw[important line] (\col,\row) coordinate (A)  -- (\col+1,\row) coordinate (B);
   				 \draw[important line] (\col,\row+1) coordinate (A)  -- (\col+1,\row+1) coordinate (B);
   				 \draw[important line] (\col,\row) coordinate (A)  -- (\col,\row+1) coordinate (B);
   				 \fill (A) circle (2pt);
   				 \fill (B) circle (2pt);
   				  \draw[important line] (\col+1,\row) coordinate (A)  -- (\col+1,\row+1) coordinate (B);
   				 \fill (A) circle (2pt);
   				 \fill (B) circle (2pt);
   		 }
   	}
	\node (Q) at (1.5,0.5) {$Q_0$};
	\node (Q) at (1.5,1.8) {$Q_1$};
	\node (Q) at (1.5,3.1) {$Q_2$};
    \draw[thick,<-] (2.2,1.3cm) -- (3.5, 1.1) node[anchor=west] {$\Omega_1$};
    \draw[thick,<-] (2.2,1.0cm) -- (3.5, 1.1);
    \draw[thick,<-] (2.2,2.3cm) -- (3.5, 2.5) node[anchor=west] {$\Omega_2$};
    \draw[thick,<-] (2.2,2.6cm) -- (3.5, 2.5);
\end{tikzpicture}}
        \subcaption{}
        \label{fig:spaceTimeDomainDPST}
    \end{subfigure}
    \caption{(a) Space-time domain with discretizations according to (b) the C-SST and (c) the D-PST approach.}
    \label{fig:spaceTimeDomainAndDiscretizations}
\end{figure}

In ST-FEM, the complete computational domain is defined via a sweeping of the computational domain along the time direction. The resulting domain is denoted as space-time domain $Q$ (cf.\ Fig.\ \ref{fig:spaceTimeDomain}). The domain $Q$ may then be discretized in two different fashions: (1) the entire domain is discretized with simplex elements (referred to as time-continuous simplex space-time (C-SST) approach and illustrated in Fig.\ \ref{fig:spaceTimeDomainCSST}), or (2) the domain is divided at individual time steps $t_n$, resulting in weakly-coupled portions, the so-called space-time slabs $Q_n$ (referred to as time-discontinuous prismatic space-time (D-PST) approach and illustrated in Fig.\ \ref{fig:spaceTimeDomainDPST}).
\par
In this work, we focus on the D-PST method. Within this approach, similar to classical time-stepping schemes, the solution is computed slab-wise. As such, the length of a slab in time direction can also be interpreted as a time-step size $\Delta t$. Following this procedure requires the inclusion of a so-called jump term (second term in Eq.\ \eqref{eq:weakFormulation}) that ensures a weak coupling of the solution between the slabs.  For ansatz and test functions $T^h$ and $v^h$ from appropriate finite-dimensional functions spaces $\mathcal{S}^h_n$ and $\mathcal{V}^h_n$, the formulation for the space-time slab $Q_n$ ultimately reads:
\bigskip\par
\textit{Find } $T^h \in \mathcal{S}^h_n$\textit{, such that }$\forall v^h \in \mathcal{V}^h_n$
\begin{equation}
    \int_{Q_n} \left( v^h \frac{\partial T^h}{\partial t} + \nabla v^h \cdot \frac{\kappa_s}{\rho_s\,c_{p,s}} \nabla T^h \right) \diff Q
    \\+\int_{\Omega_n} \left( v^h\right)^+_n \left( \left(T^h\right)^+_n-\left(T^h\right)^-_n\right) \diff\Omega
    = 0,
    \label{eq:weakFormulation}
\end{equation}
with thermal conductivity $\kappa_s$, density $\rho_s$, specific heat capacity $c_{p,s}$, and spatial domain $\Omega_n$ at time $t_n$. Note explicitly that the domain shape $\Omega_n$ will deform from one slab to the next to account for the movement of the heat source. Furthermore, we use the notation
\begin{align}
    \left( T^h \right)^{\pm}_n = \lim_{\varepsilon\rightarrow 0} T^h\left(t_n \pm \varepsilon \right),
\end{align}
to illustrate the fact that at each time instance, there exists a solution value connected to the previous space-time slab and one connected to the next space-time slab. For further details on ST-FEM, the reader is referred to \cite{Behr1992,Behr1992.2}.

\subsection{Computing the heat flux around the heat source}
\label{sec:cbfMethod}

Using the novel approach for the transient CCM process, we have shown in Section \ref{sec:transientVelocity} that the time-dependent heat flux over the contact surface $q_s$ is essential to retrieve the velocity of the heat source.

Most numerical schemes estimate boundary fluxes from the finite element solution of a problem, for instance with Gauss point evaluations of the derivatives of the basis functions at the boundary \cite{Gresho1981}. However with classical piecewise-polynomial basis functions we do not have continuity of the derivatives across elements, i.e., at nodes. Thus, we cannot directly evaluate the temperature gradient at the nodes to compute the heat flux and an additional technique to reconstruct these values is needed. We then compute the heat flux $q_s$ at each time step via the consistent boundary flux (CBF) method \cite{Zhong1993}. Starting from the space-time formulation of Eq.\ \eqref{eq:weakFormulation}, we pick different spaces of ansatz and test functions $\mathcal{S}^h_{n,\text{CBF}}$ and $\mathcal{V}^h_{n,\text{CBF}}$, such that $v^h\neq0$ on the interface of interest $\Gamma_{\text{int}}$ where we would like to evaluate the heat flux. The resulting formulation reads:
\bigskip\par
\textit{Find } $\nabla T^h\cdot\textbf{n} \in \mathcal{S}^h_{n,\text{CBF}}$\textit{, such that }$\forall v^h \in \mathcal{V}^h_{n,\text{CBF}}$ 
\begin{multline}
    \int_{Q_n} \left( v^h \frac{\partial T^h}{\partial t} + \nabla v^h \cdot \frac{\kappa_s}{\rho_s\,c_{p,s}} \nabla T^h \right) \diff Q
    +\int_{\Omega_n} \left( v^h\right)^+_n \left( \left(T^h\right)^+_n-\left(T^h\right)^-_n\right) \diff\Omega
    \\ = \int_{\Gamma_{\text{int}}}v^h\frac{\kappa_s}{\rho_s\,c_{p,s}}\nabla T^h\cdot\textbf{n}\diff\Gamma,
    \label{eq:weakCBF}
\end{multline}
where $\textbf{n}$ is the normal to the boundary $\Gamma_{\text{int}}$. Finally, we can retrieve the heat flux $q_\text{num}=-\kappa_s\nabla T^h\cdot\textbf{n}$.

Note that after solving Eq.\ \eqref{eq:weakFormulation} the left-hand side of Eq.\ \eqref{eq:weakCBF} is known. Thus the computational cost of the CBF method is negligible and it is handled as a post-processing step.

\subsection{Simulating the motion of the heat source}
\label{sec:simulationMovement}

\begin{figure}
\centering
    \includegraphics[trim={2.0cm 0cm 3.0cm 0.5cm},clip, width=0.8\textwidth]{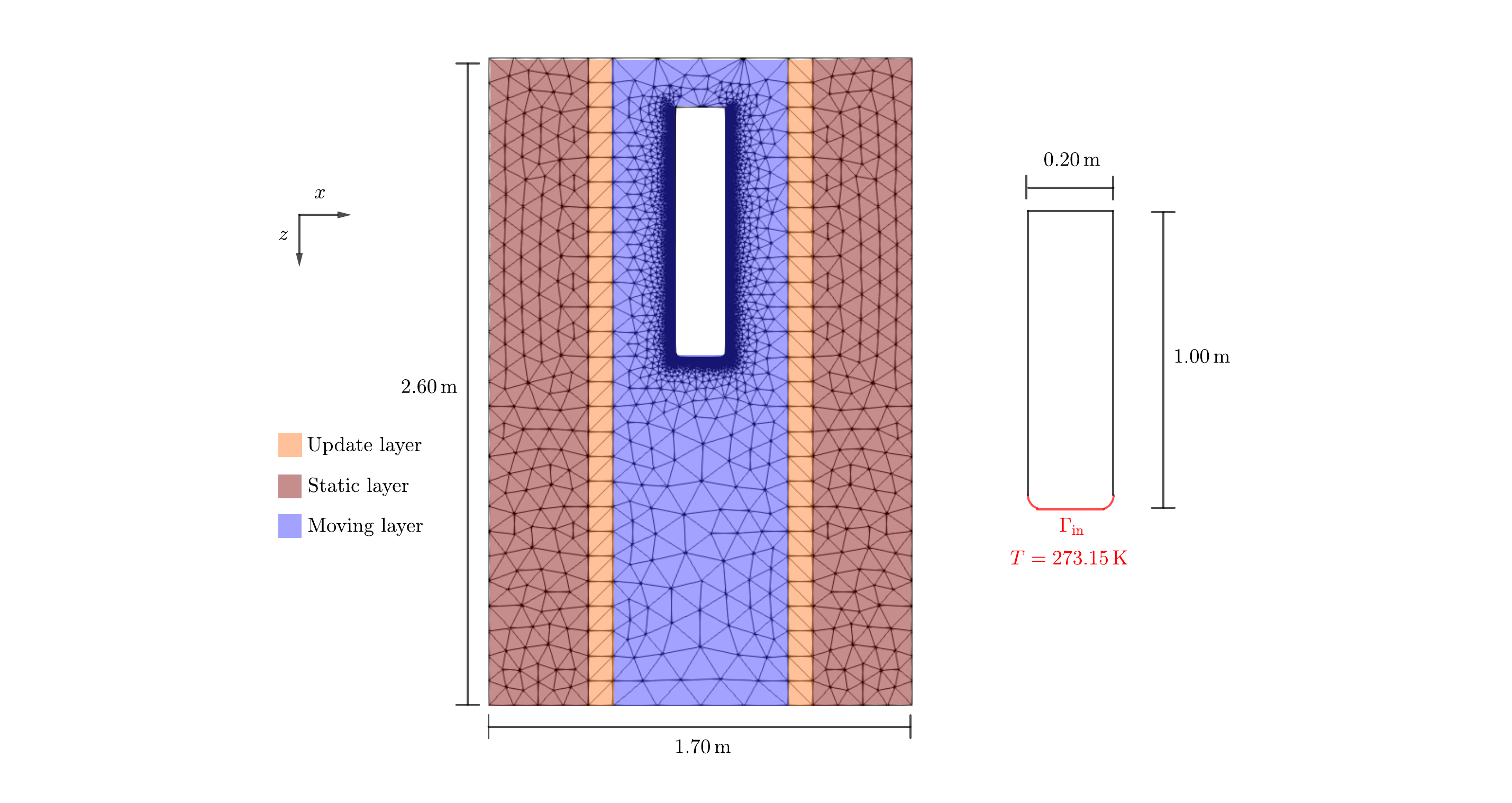}
    \caption{2D planar heat source computational domain. On the left, we show the unstructured grid and the three regions for the employed mesh-update method. The computed melting velocity is applied to the moving layer. On the right, a close-up of the heat source geometry is depicted (not drawn to scale). Note that the melting temperature is imposed on the tip of the object as shown in red.}
    \label{fig:compDomain}
\end{figure}

Given the velocity of the heat source $U$ and the time-step size $\Delta t$, we are able to compute the distance traveled by the object in a certain period of time, i.e., space-time slab. This movement has to be applied to the computational mesh in order to show the evolving position of the heat source. Since we consider the solid material to be static, the translation of the heat source results in a large relative movement with respect to the enclosing solid for longer periods of time. For such a scenario, the Virtual Region Shear-Slip Mesh Update Method (VR-SSMUM) has been introduced \cite{Key2018}. In this method, the computational mesh is split into one or multiple static and moving parts, which are each connected by a thin layer of mesh elements, the so-called update layer. In our application, the surrounding solid constitutes the static part, while the moving part is connected to the heat source instead. Fig.\ \ref{fig:compDomain} gives an example of a computational mesh with corresponding static and moving parts. Note that by employing the VR-SSMUM we can mantain the mesh refinement around the heating zone throughout the entire simulation. 
\par
What makes the mesh update method feasible is the introduction of a virtual region, i.e., an additional region of mesh elements. In particular, we connect the delimiting boundaries of the moving region, which consequently --- in an abstract sense --- form a closed ring. Note that elements in the virtual region are deactivated for the actual solution process; their purpose solely lies in modeling the mesh movement. Once the heat source starts to move downwards, mesh elements at the bottom of the computational domain drop out and, via the virtual ring, can be re-used by entering the domain from the top at a later point. This makes it possible to easily examine a large translation of the heat source, without any need for global remeshing and subsequent projection of the solution field.

\section{Verification}
\label{sec:numericalResults}

Before moving to the application of interest, we demonstrate the feasibility of the central computational techniques of our approach, namely the reconstruction of the heat flux and the effect of the moving mesh on the spatio-temporal discretization.

\subsection{Verification: heat-flux recovery}

\begin{figure}
    \centering
    \begin{subfigure}[b]{0.4\textwidth}
    \includegraphics[trim={0cm 0cm 0cm 0cm},clip, width=0.95\textwidth]{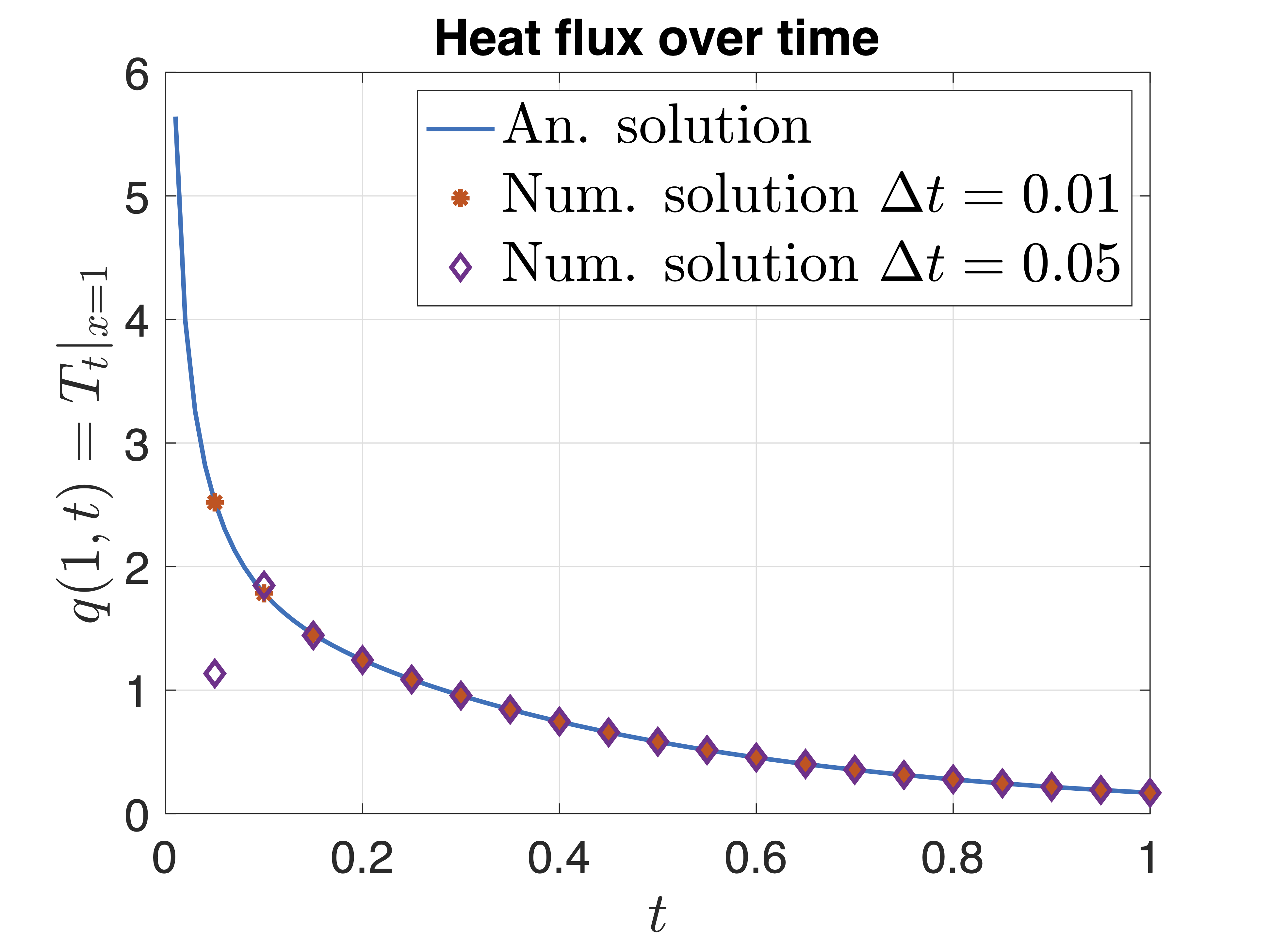}
    \caption{}
    \label{fig:cbfSolution}
    \end{subfigure}
    \begin{subfigure}[b]{0.4\textwidth}
    \includegraphics[trim={0cm 0cm 0cm 0cm},clip, width=0.95\textwidth]{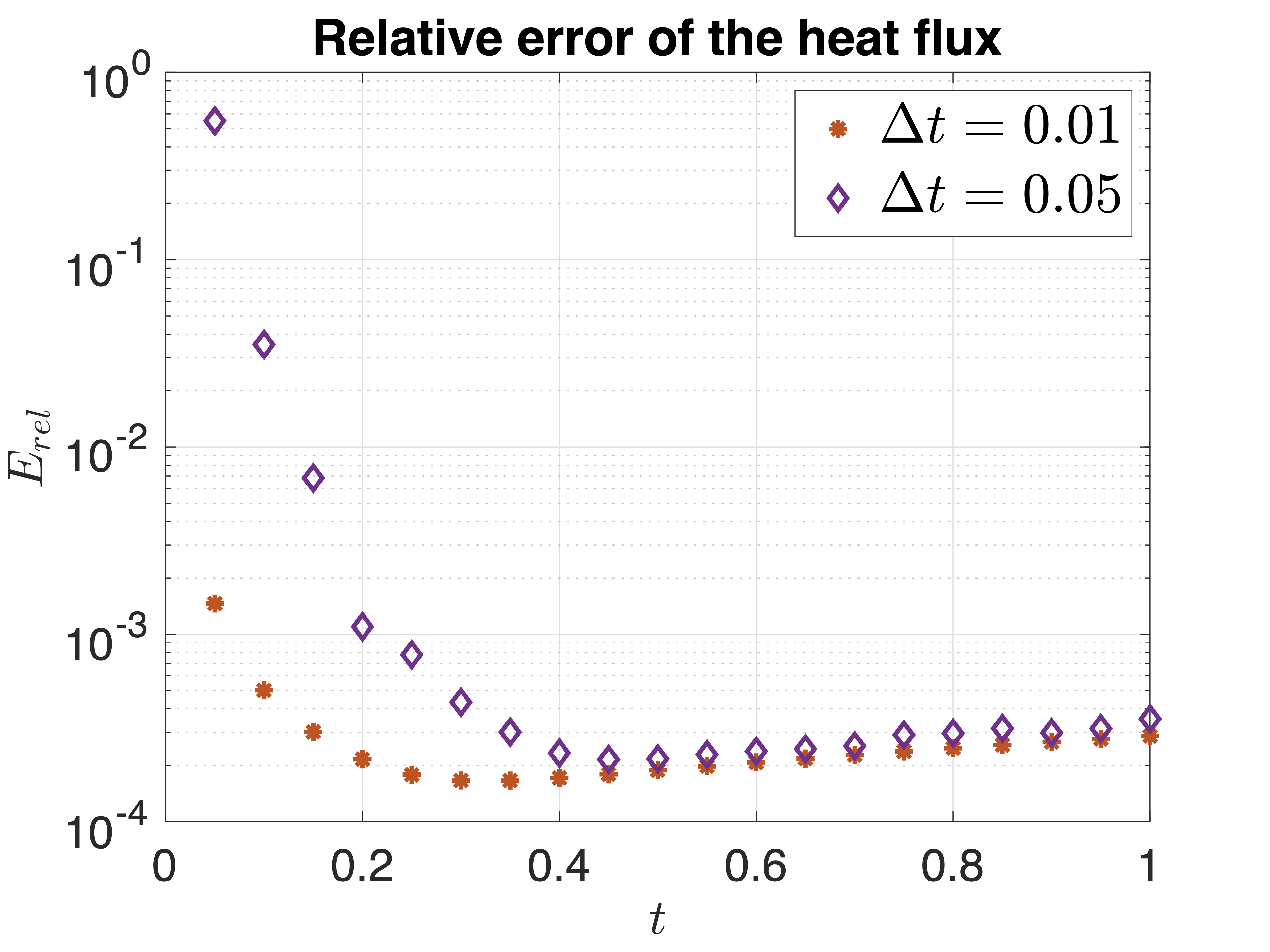}
    \caption{}
    \label{fig:cbfError}
    \end{subfigure}
    \caption{Consistent boundary flux method verification case for two different time steps. On the left, the numerical approximation of the heat flux at the right boundary is plotted against the analytical solution, see Eq.\ \eqref{eq:heatFluxAnalytical}. On the right, we show the relative error over time.}
    \label{fig:cbfHeatFlux}
\end{figure}

In a first test, we investigate the accuracy of the heat flux reconstruction at a boundary. Let us consider the heat equation over the unit square $\Omega=[0,1]^2$, such that
\begin{equation}
\begin{aligned}
    \frac{\partial T(\textbf{x},t)}{\partial t}-\Delta T(\textbf{x},t)&=0, \\ \frac{\partial T(\textbf{x},t)}{\partial x}\bigg|_{x=0} &= 0, \\ T(x=1,t)&=0, \hspace{3mm} T(x,t=0)=1.
\end{aligned}
\label{eq:verification}
\end{equation}
That is, we impose a homogeneous Neumann boundary condition on the left side of the square and a homogeneous Dirichlet condition on the right side, respectively. Note that the problem of Eq.\ \eqref{eq:verification} has a two-dimensional setting, but the corresponding temperature field depends only on the $x$-coordinate. Following \cite{mackowski2012}, we retrieve the analytical solution through separation of variables, which gives
\begin{equation}
    \hat{T}(x,t) = 2\sum_{n=1}^\infty \frac{(-1)^{n-1}\cos(\lambda_n x)}{\lambda_n x}e^{-\lambda_n^2 t}, \hspace{5mm} \lambda_n = (2n-1)\frac{\pi}{2},
\end{equation}
and the heat flux at the right boundary
\begin{equation}
    \hat{q}(x=1,t) = \frac{\partial\hat{T}}{\partial x}\bigg|_{x=1} = 2\sum_{n=1}^\infty e^{-\lambda_n^2 t}.
    \label{eq:heatFluxAnalytical}
\end{equation}
We simulate the problem on a uniform structured grid with cell size $0.02$ and we compute 20 time steps with $\Delta t=0.05$. Fig.\ \ref{fig:cbfHeatFlux} shows the computed heat flux against Eq.\ \eqref{eq:heatFluxAnalytical} and the relative error $E_{rel}(t)=|q(t)-\hat{q}(t)|/\hat{q}(t)$ over time. The results are in good agreement with the analytical solution. However, as expected, the computation of the numerical flux is problematic close to $t=0$, since the temperature gradient tends to infinity. Note that choosing a smaller time step improves the accuracy of the CBF method. In the second simulation, performed with $\Delta t=0.01$, we obtain a much lower error close to $1\times10^{-3}$ at $t=0.05$, see Fig.\ \ref{fig:cbfError}.

\subsection{Verification: mesh-update method}

\begin{figure}
    
    \centering
    \begin{subfigure}[b]{0.30\textwidth}
    \includegraphics[trim={0.5cm 0.7cm 0.5cm 0.7cm},clip, width=\textwidth]{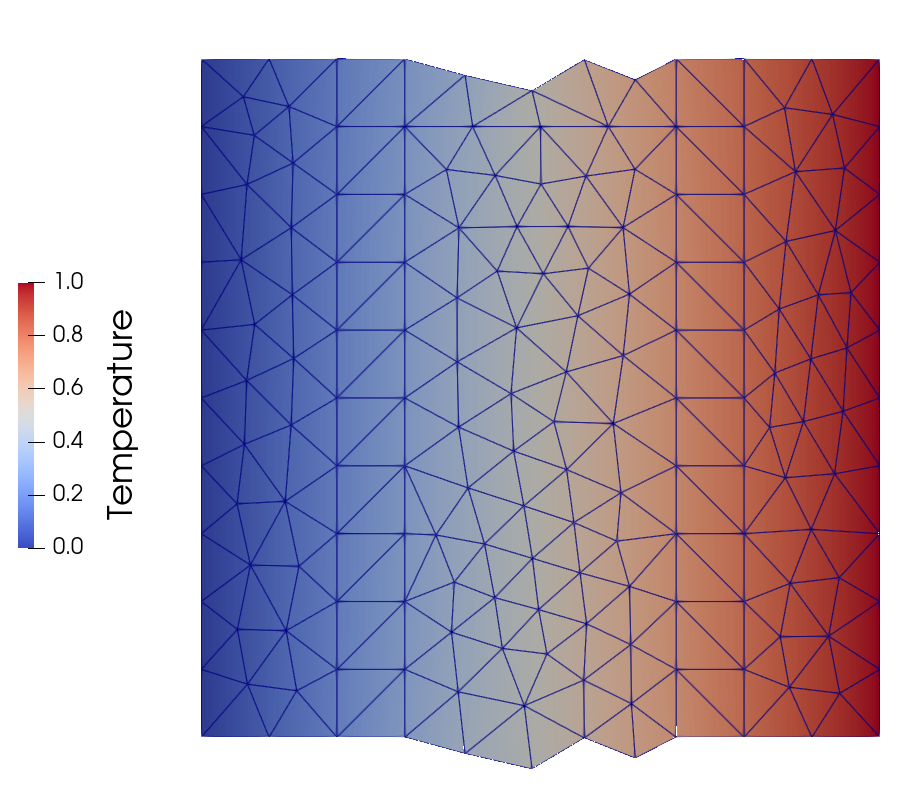}
        \subcaption{Temperature at $t=20$}
         \label{fig:vrssmumTemp}
    \end{subfigure}
    \begin{subfigure}[b]{0.30\textwidth}
    \includegraphics[trim={0.5cm 0.7cm 0.5cm 0.7cm},clip, width=\textwidth]{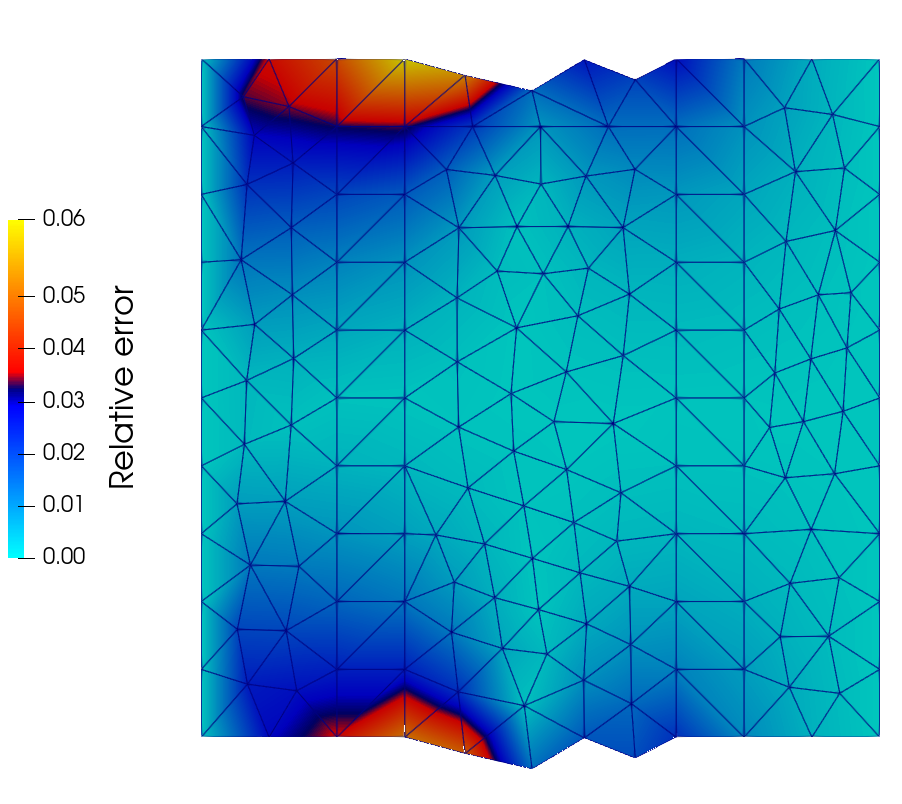}
        \subcaption{Nodal error at $t=20$}
         \label{fig:vrssmumError}
    \end{subfigure}
    \begin{subfigure}[b]{0.36\textwidth}
    \includegraphics[trim={0cm 0cm 0cm 0cm},clip, width=\textwidth]{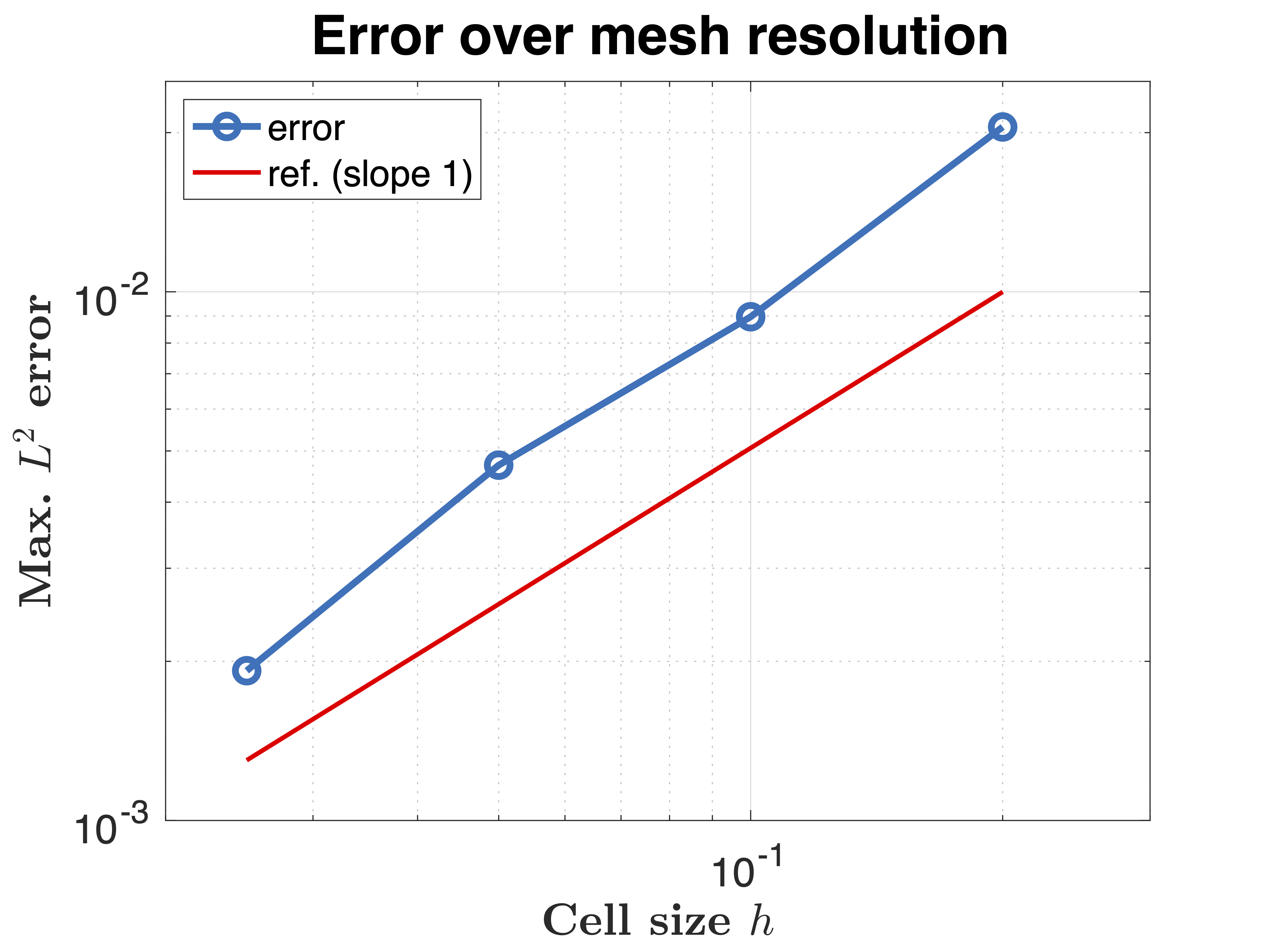}
        \subcaption{Convergence of $L^2$ error}
        \label{fig:vrssmumConvergence}
    \end{subfigure}
    
    \caption{Mesh-update method verification case. We solve for the transient heat equation on a domain where the middle part moves downwards over time. The temperature profile (a) and the corresponding nodal relative error (b), obtained on a grid with cell size $h=0.1$, are plotted at the final time step. The maximum $L^2$ error (c) is computed against four different mesh sizes.}
    \label{fig:vrssmumTestCase}
    
\end{figure}

In the second verification case, we examine to which extent the numerical solution is affected by the partially moving mesh. We consider again Eq.\ \eqref{eq:verification} on a unit square, but this time with the boundary conditions
\begin{equation}
\begin{aligned}
    \frac{\partial T(\textbf{x},t)}{\partial y}\bigg|_{y=0} &= 0, \hspace{3mm} \frac{\partial T(\textbf{x},t)}{\partial y}\bigg|_{y=1} = 0, \\ T(x=0,t)&=0, \hspace{3mm} T(x=1,t)=1, \\
    T(\textbf{x},t=0) &= x, 
\end{aligned}
\label{eq:vrssmumTestCase}
\end{equation}
so that the analytical solution is $\hat{T}(\textbf{x},t)=x$. We employ a partially moving mesh where the middle layer, i.e.\ for $0.3<x<0.7$, moves in the negative $y$-direction with velocity $0.005$. We compute 20 time steps with $\Delta t=1$ on a triangular mesh.

Fig.\ \ref{fig:vrssmumTestCase} shows the temperature profile at the final time step and the $L^2$-norm error against the analytical solution, i.e., $\max_t\lVert T-\hat{T}\rVert_{L^2(\Omega(t))}$. We consider four different unstructured meshes, where the grid cell size $h$ ranges from $0.2$ to $0.025$. As we employ a finite element discretization with piecewise linear basis functions, we should expect second order convergence for the considered heat equation, while the computed rate is close to one. Thus, it is evident that the mesh movement introduces an additional error in the elements close to the moving region boundaries, as shown in Fig.\ \ref{fig:vrssmumError}. This is caused by the activation and deactivation of the boundary elements, where we cannot enforce the weak continuity of Eq.\ \eqref{eq:weakFormulation} as there is no information on the previous time step. Note, however, that Fig.\ \ref{fig:vrssmumConvergence} shows clear convergence for refined meshes, so that the error introduced by the mesh-update method can be controlled via finer computational grids. We refer to \cite{Key2018} for additional details on the properties of the method.

\section{Application cases and numerical results}

In the previous sections, we introduced the necessary building blocks of our numerical approach to simulate transient CCM. In order to demonstrate the feasibility of our method, we consider two application cases. 

\subsection{Transient CCM applied to the sinking motion of a thermal melting probe}
\label{sec:2dProbe}

The underlying motivation for the first simulation scenario is to understand the dynamic behavior of thermal melting probes. Autonomous cryobots are a promising technology for extraterrestrial ice exploration \cite{konstantinidis2014,Kowalski2016} and are planned to be deployed on the Saturnian moon Encelauds \cite{Sherwood2016} and the polar ice caps of Mars \cite{Komle2018}. While high-fidelity models for the spatio-temporal evolution of the melt channel are investigated \cite{Boledi2022}, namely by considering the convection-coupled phase-change process around the probe, our approach allows to study the melting velocity ramp-up in response to the low temperatures of the target ice bodies.

\subsubsection{2D planar probe - displacement and temperature evolution}
\label{sec:probeDisplacement}

We now introduce the computational mesh and the general setup that we use to simulate the problem. Note that in this case the thermal melting probe, also referred to as cryobot, corresponds to the heat source, while the ambient phase-change material is the ice that the cryobot melts into. Fig.\ \ref{fig:compDomain} shows the geometry of the planar probe with a rounded melting head. The probe is placed in the moving layer, depicted in blue, where the nodal coordinates are updated over time. The remaining parts of the computational mesh are constructed to work with the mesh-update method, as described in Section \ref{sec:simulationMovement}. We emphasize that the element size is not constant, but it is smaller closer to the cryobot in order to reduce the computational error where the strongest thermal gradients are be expected. As a final remark, we recall that the mesh also features a deactivated copy of the moving layer that is not pictured for the sake of clarity.

\begin{table}
\centering
\begin{tabular}{ll|ll} 
   \textbf{Ambient ice} & & \textbf{Melting parameters} & \\ 
   \midrule
   $\rho_s$ $[\si{\kg\per\m\cubed}]$ & 921.3 & $\rho_l$ $[\si{\kg\per\m\cubed}]$ & 1000.0 \\
   $c_{p,s}$ $[\si{\joule\per\kg\per\kelvin}]$ & 1877.2 & $c_{p,l}$ $[\si{\joule\per\kg\per\kelvin}]$ & 4200.0 \\
   $\kappa_s$ $[\si{\watt\per\m\per\kelvin}]$ & 2.5428 & $\kappa_l$ $[\si{\watt\per\m\per\kelvin}]$ & 0.6 \\
   $T_s$ $[\si{\kelvin}]$ & 210.0 & $T_m$ $[\si{\kelvin}]$ & 273.0 \\
   \cmidrule{1-2}
   \textbf{Probe geometry} & & $\mu$ $[\si{\kg\per\m\per\s}]$ & 0.0013 \\
   \cmidrule{1-2}
   $R_{\text{probe}}$ $[\si{\m}]$ & 0.08 & $h_m$ $[\si{\joule\per\kg}]$ & 333700 \\
   $L_{\text{probe}}$ $[\si{\m}]$ & 1.0 & $g$ $[\si{\m\per\s\squared}]$ & 3.7 \\
   $m_{\text{probe}}$ $[\si{\kg}]$ & 25.0 & & \\
\end{tabular}
\caption{2D planar probe simulation setup. The melting conditions reproduce ice found at the polar caps of Mars. The corresponding ice properties are computed according to \cite{Ulamec2007}.}
\label{table:matProps}
\end{table}

The material properties of the ambient ice, e.g., density, are in general temperature dependent and their expressions can be found in \cite{Ulamec2007}. Following the proposition in \cite{Ulamec2007}, we select a reference temperature value $T_{\text{ref}}=0.5(T_m - T_s)$, where $T_m$ is the melting temperature and $T_s$ is the initial temperature of the ice. Thus, we evaluate the material properties at $T_{\text{ref}}$ and we perform all the computations with constant parameters. Table \ref{table:matProps} shows the configuration of the problem. We consider an environment that resembles the polar caps of Mars, where the ambient ice is at the initial temperature $T_s=\SI{210}{\kelvin}$. The geometry of the probe is inspired by the IceMole prototype described in \cite{konstantinidis2014}.

\begin{figure}
    \centering
    \begin{subfigure}[c]{0.42\textwidth}
    \includegraphics[trim={0.5cm 0cm 0.5cm 0cm},clip, width=1.1\textwidth]{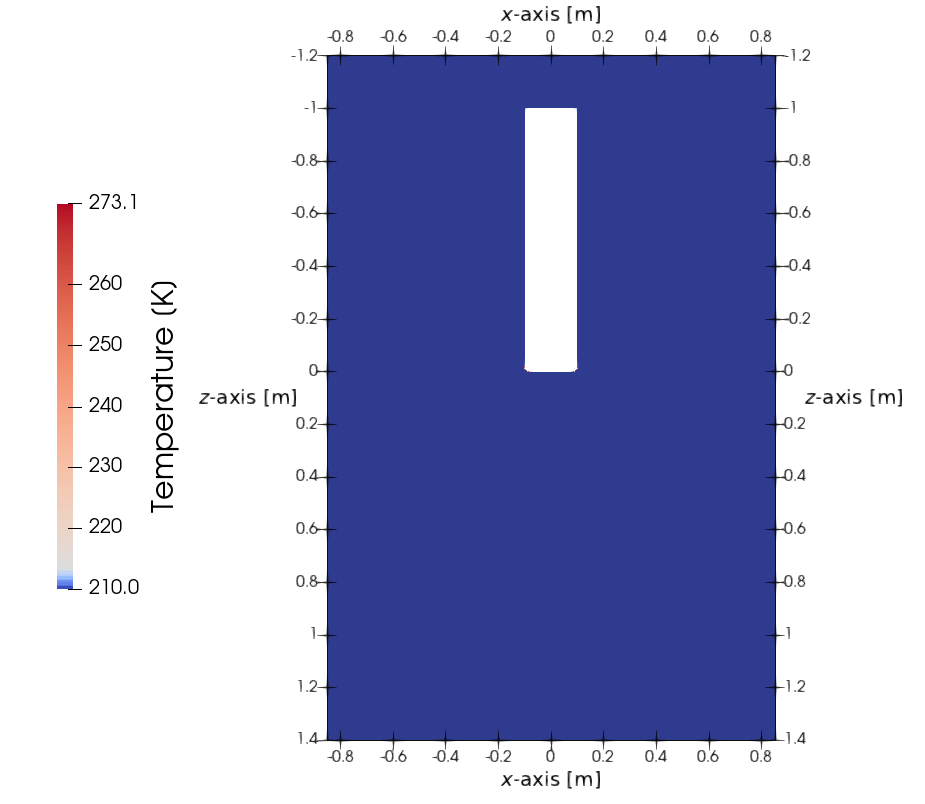}
    \caption{$t = \SI{0}{\s}$}
    \end{subfigure}
    \begin{subfigure}[c]{0.42\textwidth}
    \includegraphics[trim={0.5cm 0cm 0.5cm 0cm},clip, width=1.1\textwidth]{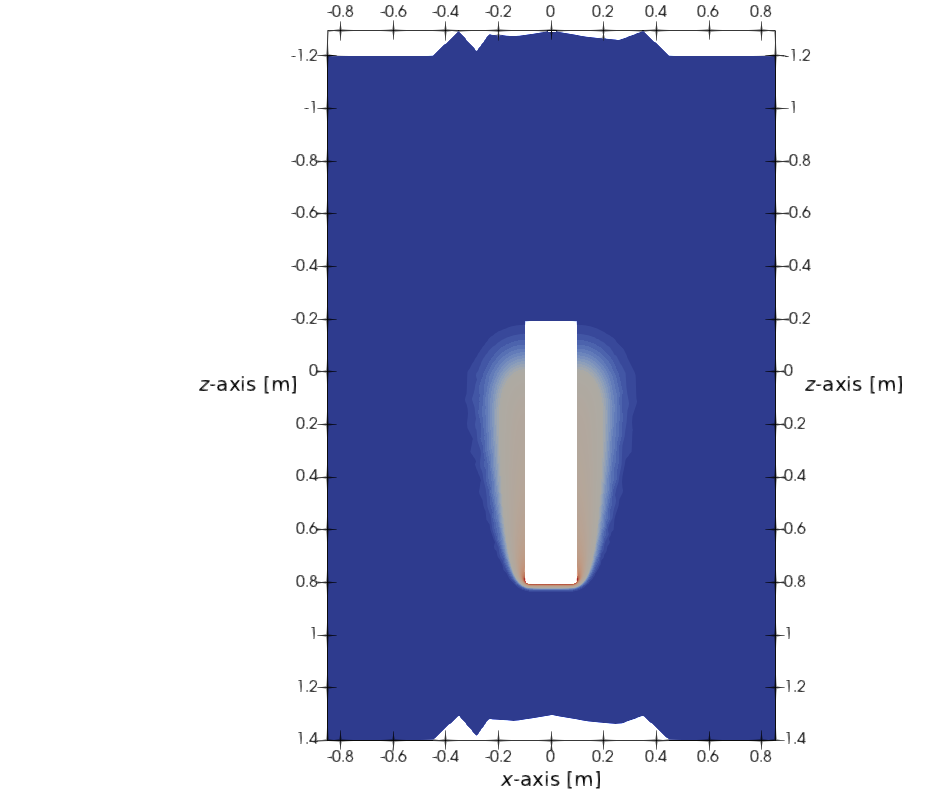}
    \caption{$t = \SI{3000}{\s}$}
    \label{fig:probeDispTempB}
    \end{subfigure}
    \caption{Temperature-controlled displacement of the probe over time due to close-contact melting. On the left, the initial state with the tip at $z=\SI{0}{\m}$ is pictured. On the right, we show the computed position and temperature distribution after $\SI{3000}{\s}$.}
    \label{fig:probeDispTemp}
\end{figure}

\begin{figure}
    \centering
    \begin{subfigure}[c]{0.42\textwidth}
    \includegraphics[trim={0.5cm 0cm 0.5cm 0cm},clip, width=1.1\textwidth]{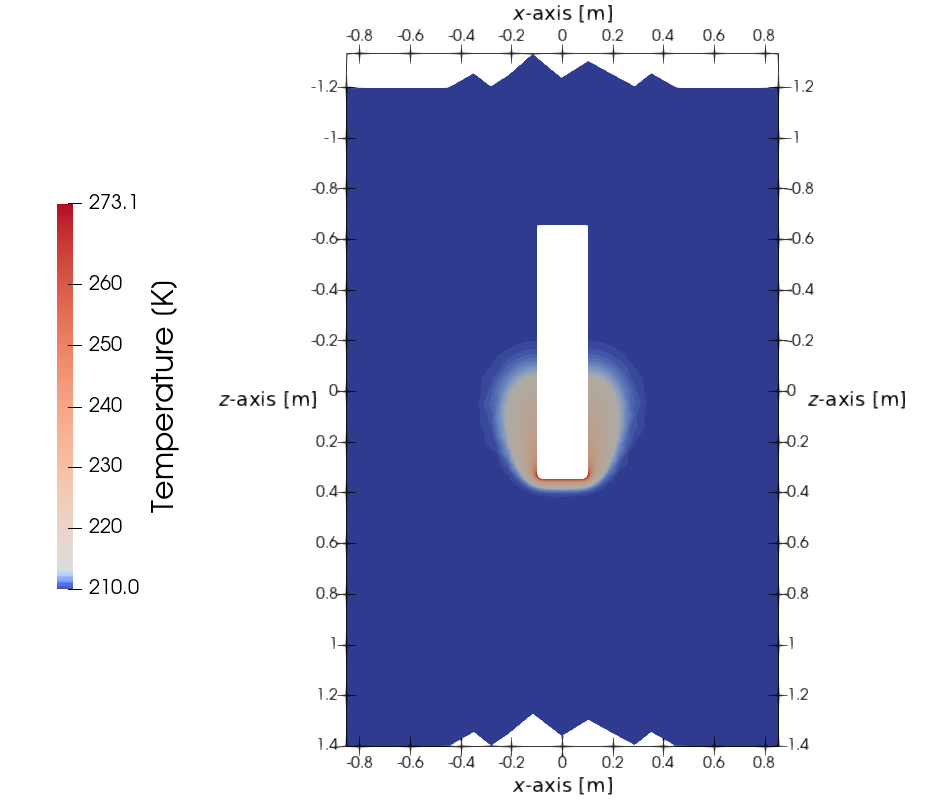}
    \caption{$\dot{Q} = \SI{1}{\kW}$}
    \label{fig:probeDispPower1}
    \end{subfigure}
    \begin{subfigure}[c]{0.42\textwidth}
    \includegraphics[trim={0.5cm 0cm 0.5cm 0cm},clip, width=1.1\textwidth]{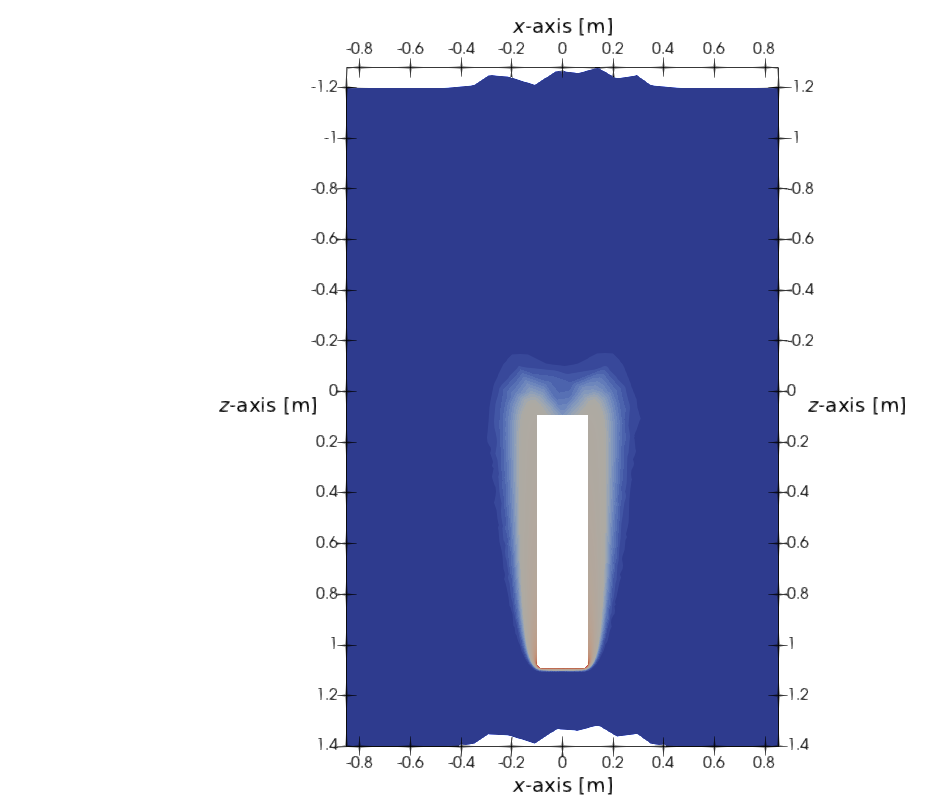}
    \caption{$\dot{Q} = \SI{4}{\kW}$}
    \label{fig:probeDispPower5}
    \end{subfigure}
    \caption{Power-controlled displacement of the probe over time due to close-contact melting. We show the computed position and temperature distribution after $\SI{3000}{\s}$ for two values of heat-flow rate.}
    \label{fig:probeDispPower}
\end{figure}

In the first simulation, we consider the equilibrium assumption of Section \ref{sec:eqVelocity}. That is, we do not yet need to recuperate the heat flux at the probe tip when computing the melting velocity. This simplifies the coupled approach presented in Section \ref{sec:transientVelocity}
and we will investigate afterwards if such simplification is valid. We consider both methods of temperature-controlled and power-controlled melting as the only difference comes from the final relation for the melting velocity. Each time, we compute 200 time steps with $\Delta t = \SI{15}{\s}$ on a mesh that comprises 60062 nodes. In order to choose the spatio-temporal resolution of the computational grid, we perform the same simulation on a finer mesh with 359116 nodes and $\Delta t=\SI{1}{\s}$, and use the result as reference. The relative error in $L^2$-norm of the temperature profile after $\SI{600}{\s}$ is equal to $4.1270\times10^{-3}$ with respect to the fine grid. We assume that the obtained numerical accuracy suffices and thus the coarser mesh is selected.

Fig.\ \ref{fig:probeDispTemp} shows the displacement of the probe over time for a temperature-controlled case, where we apply $T_w=\SI{353}{\kelvin}$ at the tip. In Fig.\ \ref{fig:probeDispTempB} the tip of the probe has reached the depth $z=\SI{0.8016}{\m}$, which is consistent with the computed melting velocity $U=\SI{2.6872e-4}{\m\per\s}$. Fig.\ \ref{fig:probeDispPower} shows the same configuration with a power-controlled computation of the melting velocity. As expected, applying a higher power at the probe tip results in a considerably faster melting process (\ref{fig:probeDispPower5}).

\begin{figure}
    \centering
    \begin{subfigure}[b]{0.35\textwidth}
    \includegraphics[trim={0cm 0cm 0cm 0cm},clip, width=0.9\textwidth]{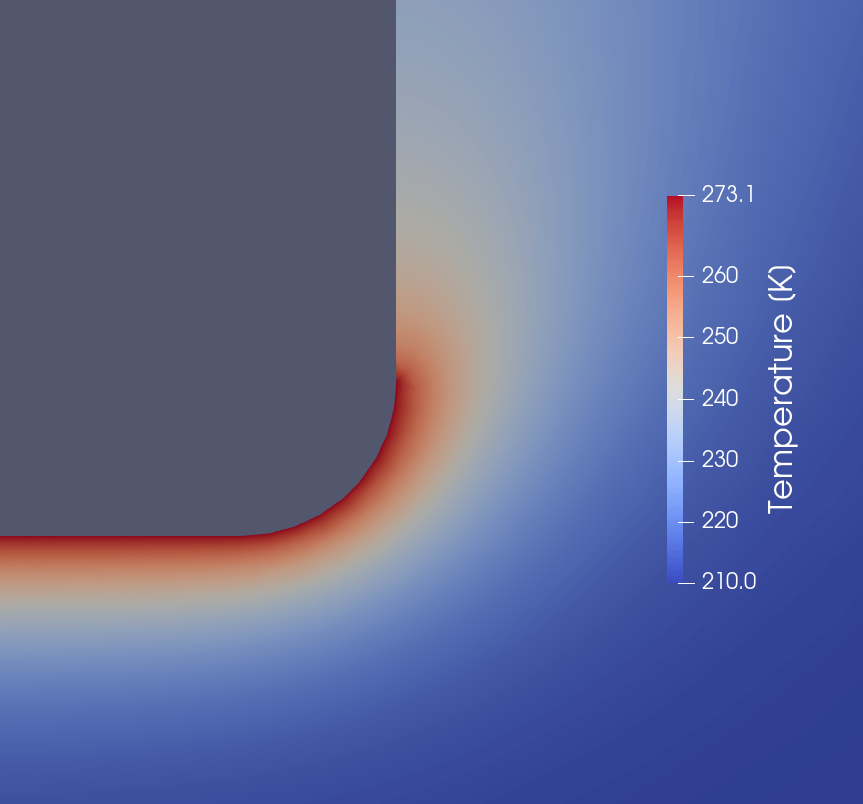}
    \caption{Close-up of Fig.\ \ref{fig:probeDispPower1}}
    \label{fig:probeTempTip}
    \end{subfigure}
    \begin{subfigure}[b]{0.4\textwidth}
    \includegraphics[trim={0cm 0cm 0cm 0cm},clip, width=1.0\textwidth]{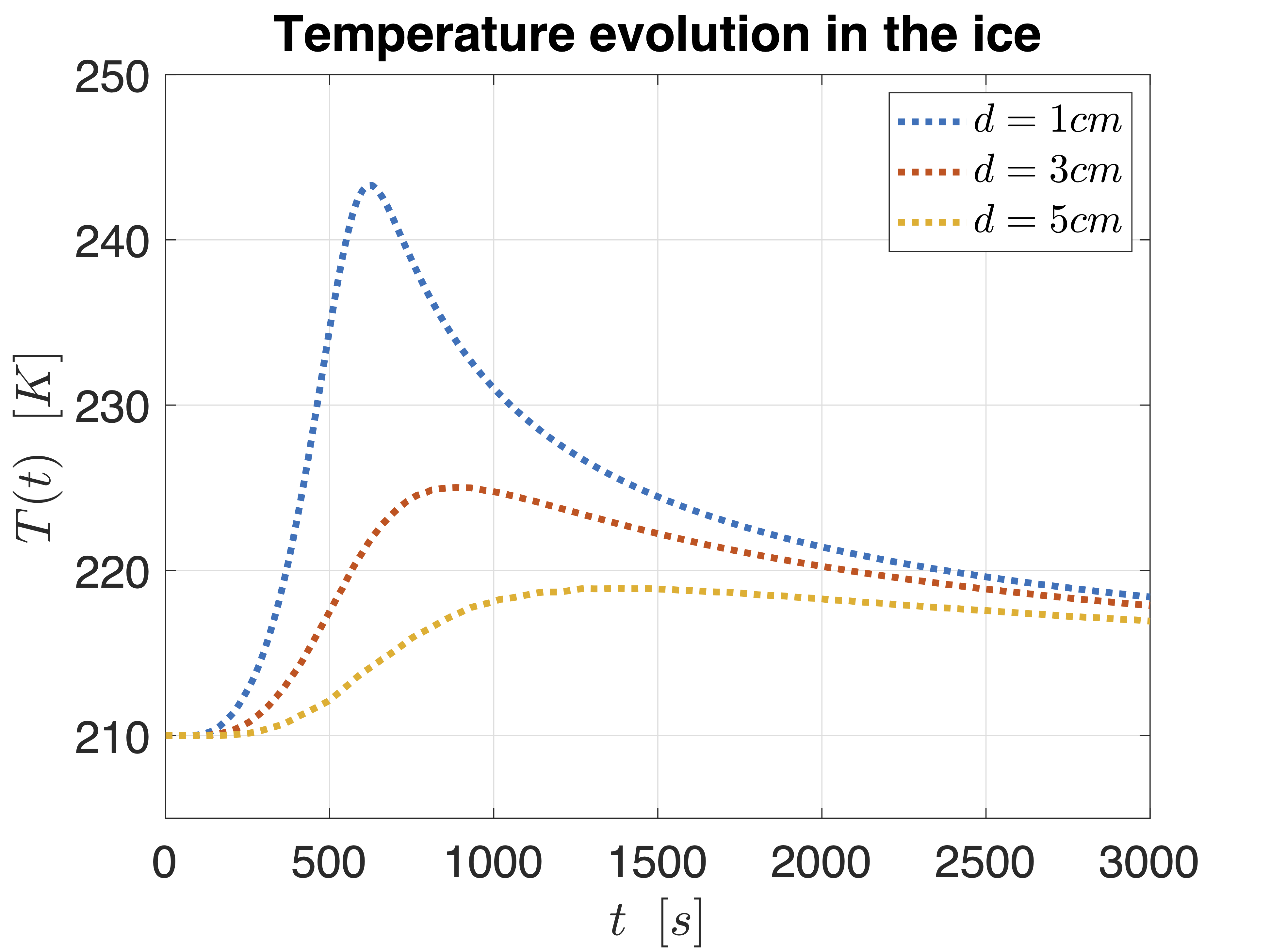}
    \caption{Temperature evolution at a distance $d$}
    \label{fig:sensorTemperature}
    \end{subfigure}
    \caption{Power-controlled temperature evolution in the ice. The transient temperature value is shown at three different locations close to the melting probe.}
    \label{fig:sensorAnalysis}
\end{figure}

We remind the reader that our approach solves the temperature equation in the surrounding ice at each time step. That is, we possess a tool for the analysis of the temperature distribution over time. We now consider a fixed point in the proximity of the probe. This setup is relevant as it replicates the data collected by a temperature sensor and could be used to validate experimental data in the future. Fig.\ \ref{fig:probeTempTip} shows a close-up of Fig.\ \ref{fig:probeDispPower1}, where the heated ice and probe tip are highlighted. Note that we employ a different color map to highlight the temperature distribution in the proximity of the probe. Then, we select three different points at $z=\SI{8}{\cm}$ and with distance $d$ equal to $\SI{1}{\cm}$, $\SI{3}{\cm}$ and $\SI{5}{\cm}$ from the probe wall. Fig.\ \ref{fig:sensorAnalysis} shows the evolving temperature at such points throughout the simulation. As expected, the temperature suddenly increases when the probe passes by the sensor and the increase is higher closer to the cryobot wall. As time progresses, the temperature tends back to the initial value $T_s$, but the cooling process is slower.

\subsubsection{2D planar probe - transient melting velocity}
\label{sec:transientTestCase}

\begin{figure}
    \centering
    \includegraphics[trim={0cm 0cm 0cm 0cm},clip, width=0.4\textwidth]{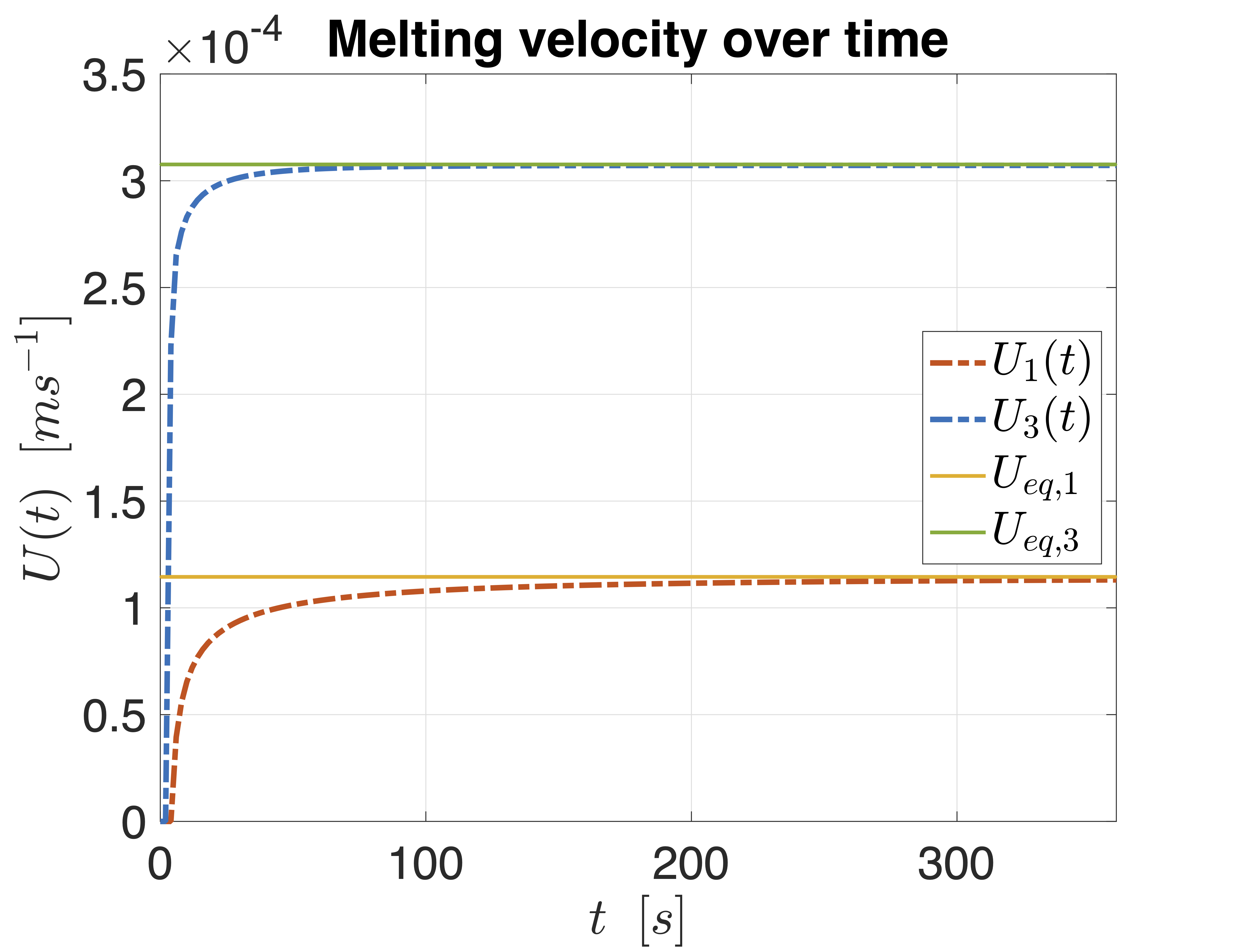}
    \caption{Evaluation of transient melting velocities for a power-controlled melting process. We consider two simulations with heat-flow rate equal to $\SI{1}{\kW}$ and $\SI{3}{\kW}$, respectively. Note that transient melting velocity asymptotically approaches the equilibrium estimates $U_{eq,i}$. The latter are computed according to Eq.\ \eqref{eq:meltPower}.}
    \label{fig:meltVelocity}
\end{figure}

Up to this point we have considered the assumption of equilibrium melting, see Section\ \ref{sec:eqVelocity}, which does not require to compute the numerical heat flux at the cryobot tip. It is important however to question if such simplification is justified. To do so, we examine the formulation presented in Section\ \ref{sec:transientVelocity}, where we utilize the numerical temperature gradient at each time step to compute the melting velocity. 

We use the same setup as in Section \ref{sec:probeDisplacement}, namely the same material properties of Table \ref{table:matProps}. The only difference comes with the considered timescale. As we are interested in the evolution of the melting velocity over time, we select a much smaller time-step size --- $\Delta t = \SI{2}{\s}$ --- and we compute 180 time steps. Before commenting on the melting velocity, the computation of the temperature gradient needs to be mentioned. After retrieving the temperature field in the ambient ice, the temperature gradient on the boundary that represents the tip of the probe is computed via the CBF method, see Section \ref{sec:cbfMethod}. That is, we now have information on the heat flux at the boundary nodes. In order to close the problem and compute the transient melting velocity, as introduced in Section \ref{sec:transientVelocity}, we simply average the numerical gradient over the number of nodes.

Fig.\ \ref{fig:meltVelocity} shows the transient melting velocities against the respective equilibrium values. At the very beginning, the computed velocity is 0, then after a rapid increase it tends to the asymptotic value $U$. This is consistent with our model, as we never exceed the equilibrium estimate of Section \ref{sec:eqVelocity}. We also point out that the two transient curves approach the asymptotic values at a different pace. We can quickly demonstrate this by considering the time needed for the melting velocity to reach 95\% of the equilibrium value. With the $\SI{1}{\kW}$ configuration, $\SI{120}{\s}$ are necessary, while only $\SI{20}{\s}$ are needed for the $\SI{3}{\kW}$ setup. 

This result entails an additional observation. The equilibrium assumption analyzed in the previous section works well on large timescales. Over the course of thousands of seconds, the error introduced in the displacement of the probe by neglecting the time-dependent velocity is negligible. However, the same cannot be said for small timescales. As shown, the melting process takes tens of seconds, depending on the configuration, to reach the equilibrium state. This workflow gives us a better understanding of such process right after starting the cryobot.

\subsection{Transient CCM applied to hot-wire cutting}
\label{sec:hotWireCutting}

Our second test case considers an industrial application, namely hot-wire cutting, where a hot wire is pressed against a low melt-point material, e.g., a polystyrene foam. Thus, the heated rod \emph{cuts} through the solid phase-change material to produce various objects, such as foam inserts for car seats \cite{Mayer2008}. Existing results in the literature employ CCM theory to investigate the process operating conditions and improve the quality of the cut.

\begin{figure}
    \centering
    \begin{subfigure}[b]{0.55\textwidth}
    \includegraphics[trim={0.5cm 2.8cm 0.5cm 3cm},clip, width=\textwidth]{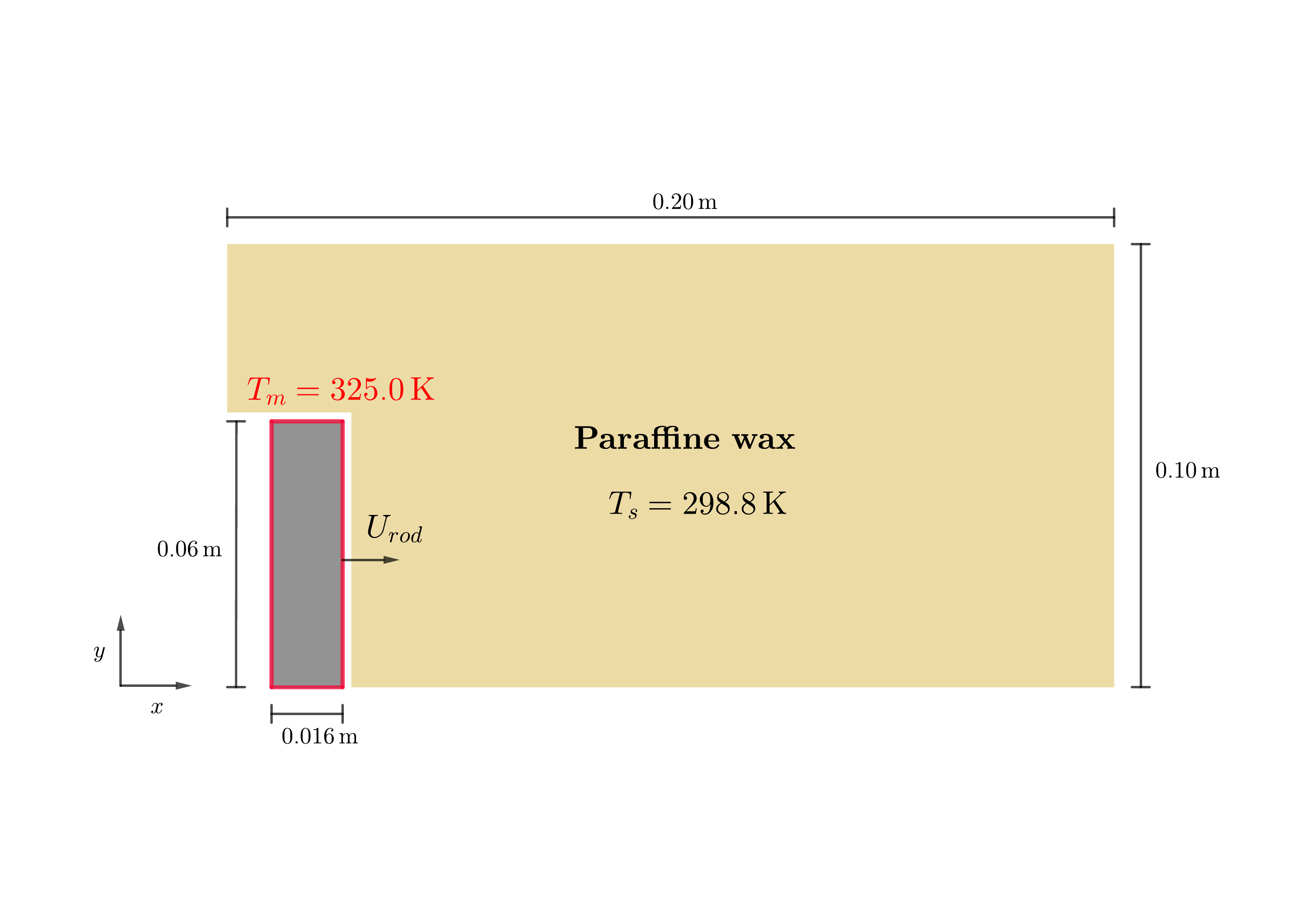}
    \caption{Setup of the simulation}
    \label{fig:waxRodSetup}
    \end{subfigure}
    \begin{subfigure}[b]{0.44\textwidth}
    \includegraphics[trim={0.5cm 0cm 0.5cm 0cm},clip, width=\textwidth]{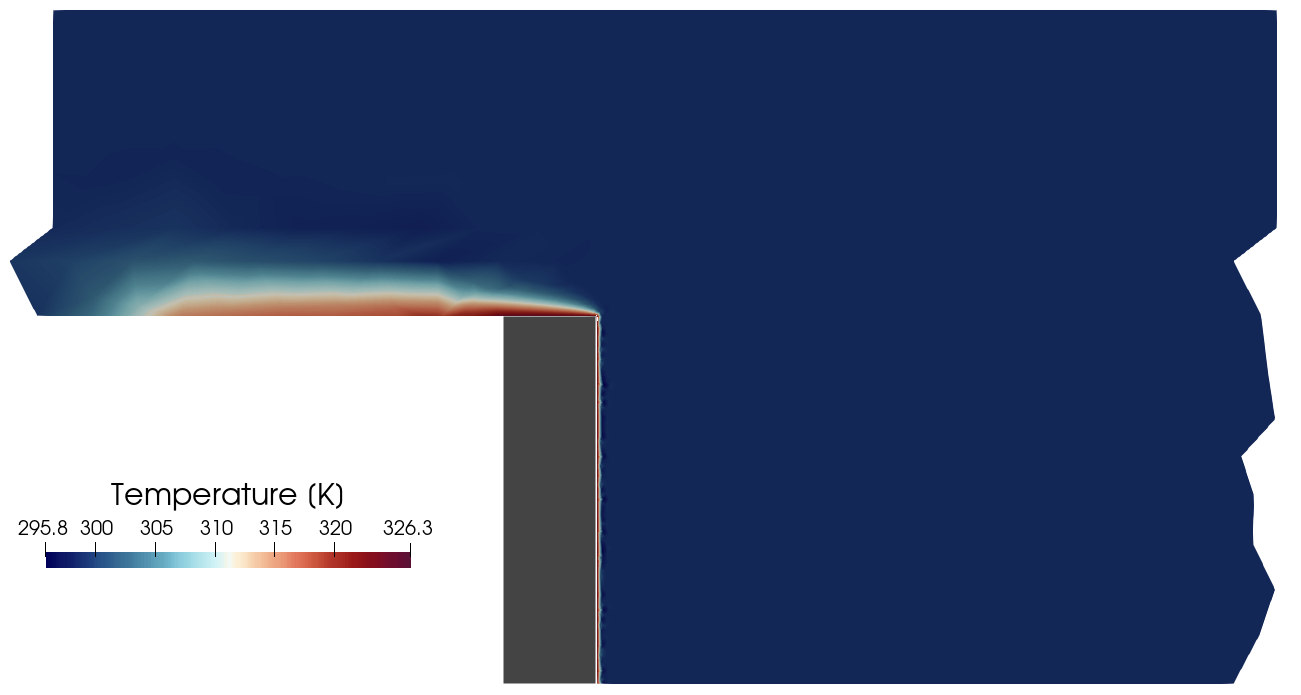}
    \caption{Rod displacement after $\SI{360}{\s}$}
    \label{fig:waxRodTemp}
    \end{subfigure}
    \caption{2D hot-wire cutting simulation. A hot copper wire is pressed against a solid block of wax, which has the initial temperature $T_s=\SI{298.8}{\K}$ (a). After $\SI{360}{\s}$ the rod has cut halfway through the block (b).}
    \label{fig:waxRodCase}
\end{figure}

\begin{table}
\centering
\begin{tabular}{ll|ll} 
   \textbf{Paraffine wax} & & \textbf{Melting parameters} & \\ 
   \midrule
   $\rho_s$ $[\si{\kg\per\m\cubed}]$ & 775.0 & $\rho_l$ $[\si{\kg\per\m\cubed}]$ & 775.0 \\
   $c_{p,s}$ $[\si{\joule\per\kg\per\kelvin}]$ & 2674.0 & $c_{p,l}$ $[\si{\joule\per\kg\per\kelvin}]$ & 2674.0 \\
   $\kappa_s$ $[\si{\watt\per\m\per\kelvin}]$ & 0.13 & $\kappa_l$ $[\si{\watt\per\m\per\kelvin}]$ & 0.13 \\
   $T_s$ $[\si{\kelvin}]$ & 298.8 & $T_m$ $[\si{\kelvin}]$ & 325.0 \\
   \cmidrule{1-2}
   \textbf{Rod geometry} & & $\mu$ $[\si{\kg\per\m\per\s}]$ & 0.00279 \\
   \cmidrule{1-2}
   $R_{\text{rod}}$ $[\si{\m}]$ & 0.03 & $h_m$ $[\si{\joule\per\kg}]$ & 221000 \\
   $L_{\text{rod}}$ $[\si{\m}]$ & 0.016 & $g$ $[\si{\m\per\s\squared}]$ & 9.81 \\
\end{tabular}
\caption{2D hot-wire cutting simulation. The material properties are inspired by \cite{Mayer2008}.}
\label{table:matPropsWaxRod}
\end{table}

The computational domain is loosely inspired by the experimental study presented in \cite{Mayer2008}. That is, we consider a $\SI{20}{\cm}\times\SI{10}{\cm}$ block of paraffine wax. A $\SI{1.6}{\cm}\times\SI{6}{\cm}$ rectangular copper rod, i.e., the heat source, is then heated to $T_{rod}=\SI{335.34}{\K}$ and pushed against the wax, which results in melting and cutting of the solid block, as shown in Fig.\ \ref{fig:waxRodSetup}. The force exerted by the rod on the wax is $F_{ex}=\SI{60}{\N}$ and we consider no gravitational contribution. Table \ref{table:matPropsWaxRod} shows the material properties for the simulation. The computational grid comprises 6582 nodes and we solve for 180 time steps with $\Delta t=\SI{2}{\s}$. Once again, the mesh is built to work with the mesh-update method, so that we have a moving layer for $y<\SI{0.07}{\m}$ where we enforce the melting velocity. To compute the latter, we consider the temperature-controlled CCM process of Eq.\ \eqref{eq:meltTempTransient}. Fig.\ \ref{fig:waxRodSetup} displays the initial setup, where the left side of the rod is located at $x=\SI{0.01}{\m}$. We recall that the melting temperature $T_m=\SI{325.0}{\K}$ is imposed on the edges of the rod as we do not solve for the liquid phase. At the end of the simulation the rod has moved $\SI{0.062}{\m}$ in the $x$-direction and the temperature field in the wax block is shown, see \ref{fig:waxRodTemp}. The computed melting velocity is $U_{rod}=\SI{1.744e-4}{\m\per\s}$.

Note that any quantitative assessment is beyond the scope of this manuscript and will be considered in future work. Nonetheless, this numerical example displays the versatility of our algorithm in handling various transient CCM processes, ranging from ice melting to industrial applications.

\section{Conclusions}
\label{sec:conclusions}

In this work we presented a transient close-contact melting model for a 2D planar heat source. We developed a numerical workflow with space-time finite elements to estimate the melting velocity of the heater and the temperature evolution in the surrounding phase-change material. The approach is based on three steps. First, we solve for the temperature field in an Eulerian reference frame around the heated body. Then, we compute the melting velocity via the numerical heat flux retrieved in the proximity of the heater. Finally, we apply this velocity to the computational mesh using the Virtual Region Shear-Slip Mesh Update Method, so that the movement of the heat source can be represented. With this technique, we distinguish between a static portion of the domain related to the solid and a moving portion connected to the heat source. The advantage is that large translations of the heater can be easily represented without the need of remeshing. The presented model results in a scale-coupled approach, as the macro-scale temperature distribution in the domain is used to compute the heat flux around the heater. The estimation of the melting velocity comes from the micro-scale close-contact melting model in the melt film. The resulting velocity is then applied to the computational mesh to represent the macro-scale movement of the heater.

After verification of the method, we studied a 2D planar thermal melting probe. This is an application-oriented simulation with a realistic geometry and material properties that resemble the ice found on the polar caps of Mars. We showed the evolving displacement of the cryobot and the spatio-temporal temperature distribution in the ambient phase-change material for equilibrium melting conditions. Furthermore, our novel scale-coupled workflow allows to directly simulate a transient close-contact melting process for the first time. In particular, we studied the velocity ramp-up after powering up the heat source and observed that the transient melting velocity increases asymptotically towards the equilibrium velocity estimate. This result enables us to quantify the needed time scale for equilibrium conditions to establish, which constitutes an improvement with respect to available models in the literature.

As a second simulation scenario, we considered an industrial application of a 2D hot-wire cutting process. We showed the displacement of a hot rod cutting through a block of paraffine wax. This setup demonstrates the versatility of our computational approach, since it can be easily applied to a different type of problem.

These results are the starting point for deeper analyses, as there is great potential in combining subscale semi-analytical models that parametrize complex contact phenomena --- here contact-melting processes --- with macro-scale temperature evolution. The flexibility and modularity of our workflow make it possible to study more complex geometries, e.g., by considering different melting heads for thermal probes, and include 3D simulations in the future.

\section*{Conflict of interest statement}

The authors declare that no commercial or financial relationships exist that could be a potential conflict of interest.

\section*{Acknowledgments}
The authors were supported by the Helmholtz Graduate School for Data Science in Life, Earth and Energy (HDS-LEE).
The work was furthermore supported by the Federal Ministry of Economic Affairs and Energy, on the basis of a decision by the German
Bundestag (50 NA 1908). The authors gratefully acknowledge the computing time granted by the JARA Vergabegremium and provided on the JARA Partition part of the supercomputer JURECA at Forschungszentrum Jülich \cite{jureca2018}.

\bibliographystyle{elsarticle-num}
\bibliography{main}

\end{document}